\documentclass[12pt,a4paper,reqno]{amsart}
\usepackage{graphicx, psfrag, amscd, amssymb}

\usepackage{booktabs}
\usepackage{mathrsfs}
\usepackage{indentfirst}
\usepackage{color}
\usepackage{graphicx,xcolor,subfigure}
\usepackage{algorithm, algorithmicx, algpseudocode}
\usepackage{pgfplots}
\pgfplotsset{compat=1.15}
\usepackage{mathrsfs}
\usepackage{mathdots}
\usetikzlibrary{arrows}
\usepackage{caption}
\usepackage{amssymb}
\usepackage{amsthm}
\usepackage{amsmath}

\newtheorem{lem}{Lemma}[section]
\newtheorem{thm}[lem]{Theorem}
\newtheorem{prop}[lem]{Proposition}
\newtheorem{cor}[lem]{Corollary}

\newcommand{\T}{{\mathcal{T}}}

\begin{document}
\title[Three New Refined Arnold Families]
{Three new Refined Arnold Families}

\author{Sen-Peng Eu}
\address{Department of Mathematics, National Taiwan Normal University, Taipei 106, Taiwan, ROC}
\email{speu@math.ntnu.edu.tw}

\author{Louis Kao}
\address{Department of Mathematics, National Taiwan Normal University, Taipei 106, Taiwan, ROC}
\email{louiskao@ntnu.edu.tw}



\thanks{S.-P. Eu is partial supported by MOST 110-2115-M-003-011-MY3 and L. Kao is partial supported by MOST 111-2811-M-003-032-MY2}

\date{\today}

\maketitle

\begin{abstract}
The Springer numbers, introduced by Arnold, are generalizations of Euler numbers in the sense of Coxeter groups. They appear as the row sums of  a double triangular array $(v_{n,k})$ of integers, $1\leq|k|\leq n$, defined recursively by a boustrophedon algorithm. We say a sequence of combinatorial objects $(X_{n,k})$ is an Arnold family if $X_{n,k}$ is counted by $v_{n,k}$. A polynomial refinement $V_{n,k}(t)$ of $v_{n,k}$, together with the combinatorial interpretations in several combinatorial structures was introduced by Eu and Fu recently.  In this paper, we provide three new Arnold families of combinatorial objects, namely the cycle-up-down permutations, the valley signed permutations and
Knuth's flip equivalences on permutations. We shall find corresponding statistics to realize the refined polynomial arrays. 
\end{abstract}
\maketitle
\section{Introduction}
 



\subsection{Springer numbers and Arnold triangle}

Let $W$ be a finite Coxeter group $W$ with set of generators $S$.
For $w\in W$, its descent set is defined by $\mathsf{Des}(w):=\{s\in S: \ell(w)<\ell(ws)\}$, where $\ell(w)$ is the length function. For each $J\subset S$, Springer~\cite{Springer_71} considered 
the set $D_J:=\{w\in W: \mathsf{Des}(w)=J\} $ and 
define $$K(W):=\max_{J\subset S} |D_J|,$$
which we called the \emph{Springer numbers} of type $W$.
It can be proved that in type $A_{n-1}$ we have $K(\mathfrak{S}_n)=E_n$, the \emph{Euler number} which counts the number of \emph{alternating permutations} $\sigma =\sigma_1\sigma_2\dots \sigma_n \in \mathfrak{S}_n$ with $\sigma_1>\sigma_2<\sigma_3>\cdots$ and therefore Springer numbers are generalizations of the Euler numbers. See~\cite{Stanley_10,EC1} for more information on Euler numbers and alternating permutations. In 1877, Seidel defined the triangular array $(E_{n,k})$ for the calculation of $E_n$ by $$E_{n,k}=E_{n,k-1}+E_{n-1,n-k+1}\quad(n\geq k\geq 2)$$ with $E_{1,1}=1, E_{n,1}=0(n\geq 2)$, and showed that $E_n=\sum_k E_{n,k}$. The refinement $E_{n,k}$ of $E_n$ is called the {\it Entringer number} since Entringer \cite{e:66} has proved that $E_{n,k}$ is the number of alternating permutations with first entry $k$. 

 The above idea of refinement can be generalized into type $B_n$ and $D_n$. For the Coxeter group $\mathfrak{B}_n$ ($\mathfrak{D}_n$, respectively) of type $B_n$ ($D_n$, respectively) Arnold~\cite{Arnold_92} proved that  the Springer numbers can be read as the row sums 
$$K(\mathfrak{B}_n)=\sum_{k\leq n}v_{n,k},\qquad \text{and} \qquad K(\mathfrak{D}_n)=\sum_{k\leq n}v_{n,-k}.$$
of the double triangle generated recursively by the boustrophedon recurrence relations
\begin{align*}
    v_{n,-k}&=v_{n,-k-1}+v_{n-1,k}\quad(n>k\geq 1),\\
    v_{n,1}&=v_{n,-1}\quad(n\geq 2),\\
    v_{n,k}&=v_{n,k-1}+v_{n-1,-k+1}\quad(n\geq k>1),
\end{align*} 
with $v_{1,1}=v_{1,-1}=1$ and $v_{n,-n}=0$ for all $n\geq 2$. Initial values are listed in Table~\ref{tab:2}.
We call $v_{n,k}$ the \emph{Arnold numbers}.

\begin{table}[!ht]
\begin{center} \begin{tabular}{c|c|ccccc|ccccc|c}
$n\backslash k$&$K(\mathfrak{D}_n)$&$-5$&$-4$&$-3$&$-2$&$-1$&1&2&3&4&5&$K(\mathfrak{B}_n)$\\ \hline
1 & 1&&&&&1&1&&&&&1\\
2 & 1&&&&0&1&1&2&&&&3\\
3 & 5&&&0&2&3&3&4&4&&&11\\
4 & 23&&0&4&8&11&11&14&16&16&&57\\
5 & 151&0&16&32&46&57&57&68&76&80&80&361\\
\end{tabular} \end{center}
\caption{The Arnold numbers $v_{n,k}$ and Springer numbers.}
\label{tab:2}
\end{table}

\subsection{Arnold family and polynomial refinements}
A {\it signed permutation} of $[n]$ is a bijection $\sigma$ of the set $[\pm n]:=-[n]\cup[n]$ onto itself such that $\sigma(-i)=-\sigma(i)$ for each $i\in[\pm n]$. For simplicity we denote $-i$ by $\bar{i}$. A signed permutation is usually denoted by its window notation $\sigma=\sigma_1\sigma_2\ldots\sigma_n$, where $\sigma_i=\sigma(i)$. Since $\mathfrak{B}_n$ and $\mathfrak{D}_n$ are respectively the permutation models of the type $B_n$ and $D_n$ Coxeter groups~\cite{Bjorner_Brenti_06}, just like $\mathfrak{S}_n$ is for type $A_{n-1}$,
by an abuse of notation we let $\mathfrak{B}_n$ denote the set of signed permutations of $[n]$, and let $\mathfrak{D}_n\subset \mathfrak{B}_n$ be  those signed permutations with $|\{i:\sigma_i<0\}|$ being even.

We may define three types of alternating signed permutations~\cite{Arnold_92,JosuatVerges_10} as
$$\mathcal{S}_n:=\{\sigma \in \mathfrak{B}_n: \sigma_1>\sigma_2<\sigma_3>\cdots \},$$
$$\mathcal{S}_n^0:=\{\sigma \in \mathfrak{B}_n: \sigma_1>0 \text{ and }\sigma_1>\sigma_2<\sigma_3>\cdots  \}.$$
$$\mathcal{S}_n^D:=\{\sigma \in \mathfrak{B}_n: \sigma_1>-\sigma_2 \text{ and }\sigma_1<\sigma_2>\sigma_3<\cdots  \}.$$
It is proved~\cite{Arnold_92, JosuatVerges_10} that
$K(\mathfrak{B_n})=|\mathcal{S}_n|$ and $K(\mathfrak{D_n})=|\mathcal{S}^D_n|=|\mathcal{S}_n|-|\mathcal{S}^0_n|.$ Arnold called those permutations in $\mathcal{S}_n^0$ and $\mathcal{S}_n^D$ the \emph{snakes} of type $B_n$ and $D_n$, respectively. Moreover, he proved that $v_{n,k}$ ($v_{n,-k}$, respectively) counts the number of snakes of type $B_n$ ($D_n$, respectively) with first entry $k$ ($-k$, respectively). 

Recently, Eu and Fu~\cite{Eu_Fu_23} gave a polynomial refinements of $v_{n,k}$ in the Arnold triangle as follows.
For $1\le |k|\le n$, define the polynomials $V_{n,k}=V_{n,k}(t)$ by 
\begin{equation*} \label{eqn:Springer-poly-recurrence}
\begin{aligned}
V_{n,-k} &= V_{n,-k-1}+t^{-1} V_{n-1,k} \quad (n>k\ge 1),\\
V_{n,1} &= t^2\, V_{n,-1} \quad (n\ge 2),\\
V_{n,k}  &= V_{n,k-1}+t\, V_{n-1,-k+1} \quad (n\ge k>1),
\end{aligned}
\end{equation*}
with $V_{1,1}=t^2$, $V_{1,-1}=1$, and $V_{n,-n}=0$ ($n\ge 2$).  Note that $v_{n,k}=V_{n,k}(1)$. 
The first few  polynomials $V_{n,k}(t)$ are listed in Table \ref{tab:Arnold-Springer-polynomials}.

\begin{table}[ht]
\caption{The Arnold--Hoffman polynomials $V_{n,k}(t)$.}
{\footnotesize
\begin{tabular}{c|ccccc}
    \hline
$n$\textbackslash $k$
        & 1 & 2 & 3 & 4 &  5  \\
    \hline
    1   &  $t^2$  &    &    &   &   \\
    2   &  $t^3$  &  $t+t^3$  &    &   &   \\
    3   &  $t^2+2t^4$  &  $2t^2+2t^4$ &  $2t^2+2t^4$ &   &  \\
    4   &  $5t^3+6t^5$  & $t+7t^3+6t^5$ & $2t+8t^3+6t^5$ &  $2t+8t^3+6t^5$ & \\
    5   &  $5t^2+28t^4+24t^6$  & $10t^2+34t^4+24t^6$ & $14t^2+38t^4+24t^6$ &  $16t^2+40t^4+24t^6$ & $16t^2+40t^4+24t^6$\\
    \hline
    \hline
$n$\textbackslash $k$
        & \multicolumn{1}{c}{$-5$} & \multicolumn{1}{c}{$-4$} & \multicolumn{1}{c}{$-3$} & \multicolumn{1}{c}{$-2$} &  \multicolumn{1}{c}{$-1$}      \\
    \hline
    1   &     &    &    &     &  1 \\
    2   &     &    &    &  0  &  $t$ \\
    3   &     &    &  0 &  $1+t^2$ & $1+2t^2$ \\
    4   &     & 0  & $2t+2t^3$ &  $4t+4t^3$ & $5t+6t^3$\\
    5   &  0  & $2+8t^2+6t^4$ & $4+16t^2+12t^4$ &  $5+23t^2+18t^4$ & $5+28t^2+24t^4$\\   
        \hline 
\end{tabular}
}
\label{tab:Arnold-Springer-polynomials}
\end{table}
We call $V_{n,k}$ the \emph{Arnold-Hoffman polynomials} since it is proved in \cite{Eu_Fu_23} that
$$Q_n(t)=\frac{1}{t}(V_{n,1}+V_{n,2}+\cdots+V_{n,n})$$ and 
$$P_n(t)-tQ_n(t)=V_{n,-1}+V_{n,-2}+\cdots+V_{n,-n}$$ for $1\leq k\leq n$, where the polynomials $P_n,Q_n$ defined by
$$\frac{d^n}{dx^n}\tan(x)=P_n(\tan(x))\quad \mbox{and}\quad \frac{d^n}{dx^n}\sec(x)=Q_n(\tan(x))\sec(x)$$
are introduced by Hoffman~\cite{Hoffman_99}.

Following~\cite{Shin_Zeng_21}, a sequence of combinatorial objects $(X_{n,k})$ is called an {\it Arnold family} if  $|X_{n,k}|= v_{n,k}$ for $1\leq|k|\leq n$. Furthermore, it is a \emph{refined Arnold family}
if we can find a statistic over $(X_{n,k})$ such that $V_{n,k}$
is the enumerator with respect to this statistic.
For example, for signed permutations, it is proved by 
Josuat-Verg\'es~\cite{JosuatVerges_10} that
$V_{n,k}(t)$ for $k>0$ ($k<0$, respectively) counts the type $B_n$ snakes ($D_n$, respectively) with respect to the statistic `change of sign'. Recently in~\cite{Eu_Fu_23} 
Eu and Fu came up with several new refined Arnold families, among them the signed increasing 1-2 trees, signed Andr\'{e} permutations, signed Simsun permutations of type I and type II, and complete increasing binary trees.


\section{Structures investigated and Main Theorems}

Our motivation comes from the fact that there are still other combinatorial structures counted by Euler numbers or Springer numbers, but their complete profiles (or even corresponding definitions in types $B_n$ or $D_n$) in the sense of refined Arnold families are missing. 
In this paper we investigate three of them, namely the cycle-up-down permutations, the valley signed permutations and Knuth's flip equivalence classes on permutations, and complete their profiles. \\

\subsection{Cycle-up-down permutations}~\label{CUDBD}
A permutation $\sigma \in \mathfrak{S}_n$ can be written in the \emph{standard cycle notation}, which means that in each cycle the first entry is the smallest entry of this cycle, and these smallest first entries among the cycles are increasing. For example, $\sigma=(1,8,5)(2,4)(3)(6,9,7)$ is a standard cycle notation. From now on, when we mention a cycle notation it is always standard. A permutation $\sigma\in \mathfrak{S}_n$ is called  \emph{cycle-up-down} if, in its standard cycle notation, each cycle $(\sigma_{i_1}, \sigma_{i_2}, \dots , \sigma_{i_r})$ is reverse alternating, i.e., $$\sigma_{i_1}< \sigma_{i_2}> \sigma_{i_3} < \cdots .$$ 
Denote $\mathcal{CUD}_n$ the set of cycle-up-down permutations of length $n$. They are first defined and investigated by Deutsch and Elizalde~\cite{Elizalde_Deutsch_11}. 

A signed permutation can also be written in the cycle notation
in a clear way. For example, 
$\sigma'= \bar{2}\bar{4}31\bar{6}\bar{7}5$ has the cycle notation $(1,\bar{2}, 4)(\bar{1},2, \bar{4})(3)(\bar{3})(5,\bar{6}, 7)(\bar{5},6,\bar{7})$ and
 $\sigma''=\bar{2} 4 \bar{5} \bar{1} 9 6 3 \bar{8} 7$ has the cycle notation 
$(1, \bar{2}, \bar{4})(\bar{1},2,4)(3, \bar{5}, \bar{9}, \bar{7}, \bar{3}, 5, 9, 7)(6)(\bar{6})(8,\bar{8})$. Sometimes people use the notations $((3,\bar{5},\bar{9},\bar{7})):=(3, \bar{5}, \bar{9}, \bar{7})(\bar{3}, 5, 9, 7)$ and $[3,\bar{5},\bar{9},\bar{7}]:=(3, \bar{5}, \bar{9}, \bar{7}, \bar{3}, 5, 9, 7)$ to distinguish two types of cycles in signed permutations. In these notations of cycle types, we call a signed permutation \emph{special} if its cycle notation contains no bracket cycles. Hence $\sigma'=((1,\bar{2}, 4))((3)((5,\bar{6}, 7))$ is special and $\sigma''=((1, \bar{2}, \bar{4}))[3, \bar{5}, \bar{9}, \bar{7}]((6))[8]$ is not. However, in order to simplify the representations, we just use $(a_{1,1},a_{1,2},\ldots)$ to represent $((a_{1,1},a_{1,2},\ldots))$ in the rest of this paper. That is, $\sigma'$ and $\sigma''$ will be written in $(1,\bar{2}, 4)(3)(5,\bar{6}, 7)$ and $(1, \bar{2}, \bar{4})(3, \bar{5}, \bar{9}, \bar{7}, \bar{3}, 5, 9, 7)(6)(8,\bar{8})$, respectively.

For a special signed permutation $\sigma$ we let $|\sigma|$ be the cycle notation (of a possibly new permutation) obtained by taking absolute values for every entry in the cycle notation of $\sigma$. For example $|\sigma'|=(1,2,4)(3)(5,6,7)$.

To our ends we define types $B_n$ and $D_n$ counterparts of the cycle-up-down permutations. 
A signed permutation $\sigma\in \mathfrak{B}_n$ is a  \emph{signed cycle-up-down permutation} of type $B_n$
if it simultaneously meets the following conditions.
\begin{enumerate}
\item it is special, 
\item in its cycle notation 
$\sigma=(a_{1,1},a_{1,2},\ldots)(a_{2,1},a_{2,2},\ldots)\ldots(a_{m,1},a_{m,2},\ldots)$ we have $a_{i,1}>0$ for all $i$,
\item $|\sigma| \in \mathcal{CUD}_n$. 
\end{enumerate}
Let $\mathcal{CUD}^{(B)}_n\subset \mathfrak{B}_n$ be the set of signed cycle-up-down permutations of type $B_n$. Also
define 
$$\mathcal{CUD}^{(B)}_{n,k}:=\{\sigma\in\mathcal{CUD}^{(B)}_n,a_{m,1}=k\}$$ for $1\leq k\leq n$. For example
\begin{align*}
\mathcal{CUD}^{(B)}_1=&\{(1)\};\\
\mathcal{CUD}^{(B)}_2=&\{(1,2),(1,\bar{2}),(1)(2)\};\\
\mathcal{CUD}^{(B)}_3=&\{(1,3,2),(1,3,\bar{2}),(1,\bar{3},2),(1,\bar{3},\bar{2}),(1,2)(3),(1,\bar{2})(3),\\
~& ~(1,3)(2),(1,\bar{3})(2),(1)(2,3),(1)(2,\bar{3}),(1)(2)(3)\}.
\end{align*}

A signed permutation $\sigma\in \mathfrak{B}_n$ is a  \emph{signed cycle-up-down permutation} of type $D_n$
if its cycle notation has the form 
$\sigma=(a_{1,1},a_{1,2},\ldots)(a_{2,1},a_{2,2},\ldots)\ldots(a_{m,1},\overline{a_{m,1}})$
and satisfies simultaneously:
\begin{enumerate} 
\item all cycles are not of the form $(a,\bar{a})$ except $(a_{m,1},\overline{a_{m,1}})$,
\item $a_{i,1}>0$ for all $i$, and 
\item $|\sigma|$ is a cycle-up-down permutation (with the convention $(|a_{m,1}|,|\overline{a_{m,1}}|):=(a_{m,1})$).
\end{enumerate}
Let $\mathcal{CUD}^{(D)}_n\subset \mathfrak{B}_n$ be the set of signed cycle-up-down permutations of type $D_n$
and define $$\mathcal{CUD}^{(D)}_{n,k}:=\{\sigma\in\mathcal{CUD}^{(D)}_n,a_{m,1}=k\}.$$ for $1\leq k\leq n$.
For example
\begin{align*}
\mathcal{CUD}^{(D)}_1=&\{(1,\bar{1})\};\\
\mathcal{CUD}^{(D)}_2=&\{(1)(2,\bar{2})\};\\
\mathcal{CUD}^{(D)}_3=&\{(1,3)(2,\bar{2}),(1,\bar{3})(2,\bar{2}),(1,2)(3,\bar{3}),(1,\bar{2})(3,\bar{3}),(1)(2)(3,\bar{3})\}.
\end{align*}

For $\sigma\in\mathcal{CUD}^{(B)}_n$ or $\mathcal{CUD}^{(D)}_n$ in its cycle notation, we define the
statistic {\it negative-peak} $\mathsf{npk}$ by
$$\mathsf{npk}(\sigma):=\#\{a_{i,j}:a_{i,j}<0, |a_{i,j}| \ge |a_{i,j-1}|\},$$
with the convention that for $\sigma\in \mathcal{CUD}^{(D)}_n$ the $i$ in the cycle unique $(i,\bar{i})$ also contributes $1$ to $\mathsf{npk}$. For example, 
for $\sigma'= (1,\bar{3},\bar{2})(4)(5,\bar{6})(7,9,\bar{8})\in \mathcal{CUD}^{(B)}_9$ we have
$\mathsf{npk}(\sigma')=2$ (namely $\bar{3}$ and $\bar{6}$), while
for $\sigma''= (1,\bar{9},\bar{2})(3,4)(5,8,\bar{6})(7,\bar{7})\in \mathcal{CUD}^{(D)}_9$
 we have $\mathsf{npk}(\sigma')=2$ (namely $\bar{9}$ and $\bar{7}$). Our first main result gives a realization of $V_{n,k}(t)$ in terms of signed cycle-up-down permutations.

\begin{thm}\label{cudpoly}(\textbf{Main Theorem I})
For $1\leq k\leq n$, we have
$$V_{n,k}(t)=\sum_{\sigma\in\mathcal{CUD}^{(B)}_{n,n-k+1}}t^{n+1-2\cdot \mathsf{npk}(\sigma)}$$ and $$V_{n,-k}(t)=\sum_{\sigma\in\mathcal{CUD}^{(D)}_{n,n-k+1}}t^{n+1-2\cdot \mathsf{npk}(\sigma)}.$$
\end{thm}

In other words, $$(\mathcal{CUD}^{(D)}_{n,1},\,\ldots, \,\mathcal{CUD}^{(D)}_{n,n},\, \mathcal{CUD}^{(B)}_{n,n}, \dots, \,\mathcal{CUD}^{(B)}_{n,1})$$ together with the statistic  $\mathsf{npk}$ is a realization of the refined Arnold family. The proof will be put in Section 3.

\subsection{Valley signed permutations}\label{VALBD}
Valley signed permutations are first defined by Josuat-Verg\'es, Novelli, and Thibon~\cite{JosuatVerges_Novilli_Thibon_12} and are only defined in types $B_n$ and $D_n$.

For $\sigma\in\mathfrak{S}_n$, an entry $\sigma_i$ is called a {\it valley} of $\sigma$ if $\sigma_{i-1}>\sigma_i<\sigma_{i+1}$ for $i=2,\ldots,n-1$; a {\it peak} of $\sigma$ if $\sigma_{i-1}<\sigma_i>\sigma_{i+1}$ for $i=2,\ldots,n-1$. We also call $\sigma_1$ a valley if $\sigma_1<\sigma_2$ and $\sigma_n$ a peak if $\sigma_n>\sigma_{n-1}$.
A signed permutation $\sigma=\sigma_1\sigma_2\ldots\sigma_n \in \mathfrak{B}_n$ is a \emph{valley signed permutation} of type $B_n$ if $\sigma_i<0$ implies that
 $|\sigma_{i-1}|$ is a valley in $|\sigma|$. Let 
$\mathcal{VS}^{(B)}_n$ be the set of  valley signed permutations of type $B_n$ and  
$$\mathcal{VS}^{(B)}_{n,k}:=\{\sigma\in\mathcal{VS}^{(B)}_n,\sigma_1=k\}$$ for $1\leq k\leq n$.
For example, 
\begin{align*}
 \mathcal{VS}^{(B)}_1=&\{1\}; \\
  \mathcal{VS}^{(B)}_2=&\{12, 1\bar{2}, 21\};\\
    \mathcal{VS}^{(B)}_3=&\{123, 1\bar{2}3, 132, 1\bar{3}2,213, 21\bar{3}, 231, 2\bar{3}1, 312, 31\bar{2}, 321\}.
\end{align*}

A signed permutation $\sigma=\sigma_1\sigma_2\ldots \sigma_n \in \mathfrak{B}_n$ is a \emph{valley signed permutation} of type $D_n$ if
$\sigma_1<0$, $|\sigma_1|>\sigma_2>0$ and, for $i\geq 3$, $\sigma_i<0$ implies that $|\sigma_{i-1}|$ is a valley in $|\sigma|$. Let $\mathcal{VS}^{(D)}_n$ be the set of valley signed permutations of type $D_n$ and 
$$\mathcal{VS}^{(D)}_{n,k}:=\{\sigma\in\mathcal{VS}^{(D)}_n,\sigma_1=\bar{k}\}$$ for $1\leq k\leq n$. 
We have
\begin{align*}
    \mathcal{VS}^{(D)}_1=&\{\bar{1}\};\\
    \mathcal{VS}^{(D)}_2=&\{\bar{1}2\};\\
    \mathcal{VS}^{(D)}_3=&\{\bar{2}13, \bar{2}1\bar{3},\bar{3}12, \bar{3}1\bar{2}, \bar{3}21\}.
\end{align*}

Define the statistic $\mathsf{neg}$ of $\sigma\in \mathcal{VS}^{(B)}_n$ (or $\mathcal{VS}^{(D)}_n$) by 
$$\mathsf{neg}(\sigma):=\#\{i:\sigma_i<0\}$$ and we obtain another realization of $V_{n,k}(t)$.

\begin{thm}\label{vspoly}(\textbf{Main Theorem II})
For $1\leq k\leq n$, we have
$$V_{n,k}(t)=\sum_{\sigma\in\mathcal{VS}^{(B)}_{n,n-k+1}}t^{n+1-2\cdot \mathsf{neg}(\sigma)}$$
 and 
$$V_{n,-k}(t)=\sum_{\sigma\in\mathcal{VS}^{(D)}_{n,n-k+1}}t^{n+1-2\cdot \mathsf{neg}(\sigma)}.$$
\end{thm}

In other words, $$(\mathcal{VS}^{(D)}_{n,1},\ldots,\mathcal{VS}^{(D)}_{n,n},\mathcal{VS}^{(B)}_{n,n},\ldots,\mathcal{VS}^{(B)}_{n,1})$$ together with the statistic  $\mathsf{neg}$, is a realization of the refined Arnold family. The proof will be put in Section 4.

\subsection{Knuth's flip equivalence classes}\label{KNUBD}
Let $\sigma\in\mathfrak{S}_n$ be a permutation with one-line notation $\sigma_1\ldots\sigma_k\sigma_{k+1}\ldots\sigma_n$. 
A permutation $\sigma_k\ldots\sigma_1\sigma_{k+1}\ldots\sigma_n$ is said to be a {\it flip} from $\sigma$ if 
$$k=n \quad \text{ or } \quad \sigma_{k+1}<\min\{\sigma_i:1\leq i\leq k\}.$$ 
Two permutations are {\it flip equivalent} if they can be obtained from each other by a sequence of flips. For example, in $\mathfrak{S}_3$ there are two flip equivalence classes, namely $\{123, 321, 231, 132\}$ and $\{213, 312\}$.
The following elegant result was by Knuth~\cite{Knuth_09}. 
\begin{prop}\label{flipSn}
The number of flip equivalence classes is the Euler number $E_n$.
\end{prop}

We extend the concept of flip equivalent to other types. For a signed permutation $\sigma=\sigma_1\ldots\sigma_k\sigma_{k+1}\ldots \sigma_n$ in its window notation, the signed permutation $\sigma_k\ldots \sigma_1\sigma_{k+1}\ldots \sigma_n$ is said to be a {\it flip} from $\sigma$ if 
$$k=n \quad \text{ or } \quad |\sigma_{k+1}|<\min\{|\sigma_i|:1\leq i\leq k\}.$$ 
Two signed permutations are {\it flip equivalent} if they can be obtained from each other by a sequence of flips. Denote by $[\sigma]$ the equivalence class contains $\sigma$. For example, in $\mathfrak{B}_3$, there are $16$ flip equivalence classes: $[123],[213],[\bar{1}23],[2\bar{1}3],[1\bar{2}3],[\bar{2}13],[12\bar{3}],[21\bar{3}],[\bar{1}\bar{2}3],$ $[\bar{2}\bar{1}3],[\bar{1}2\bar{3}],[2\bar{1}\bar{3}],[1\bar{2}\bar{3}],[\bar{2}1\bar{3}],[\bar{1}\bar{2}\bar{3}],[\bar{2}\bar{1}\bar{3}]$.

Note that from this definition the numbers of flip equivalences classes on signed permutations are not Springer numbers, but certain subsets are if we adapt the following. For a signed word $\sigma=\sigma_1\sigma_2\ldots\sigma_m$ such that $|\sigma_i|\neq|\sigma_j|$ for all $1\leq i,j\leq m$, define the statistic {\it signed maximum} $\mathsf{smax}(\sigma)$ of $\sigma$ by writing $\sigma$ as $\sigma=\sigma_L\hat{\sigma} \sigma_R$, where $|\hat{\sigma}|=\min(|\sigma|)$
and computing $\mathsf{smax}(\sigma)$ recursively:
\begin{enumerate}
\item If $\sigma_L=\sigma_R=\emptyset$, then $\mathsf{smax}(\sigma)=\hat{\sigma}$.
\item If $\sigma_L=\emptyset$ and $\sigma_R\not=\emptyset$, then

\begin{enumerate}
\item[---]  if $(\hat\sigma)>0$, then $\mathsf{smax}(\sigma)=\hat{\sigma}$.
\item[---]   if $(\hat\sigma)<0$, then $\mathsf{smax}(\sigma)=\mathsf{smax}(\sigma_R)$.
\end{enumerate}

\item If $\sigma_L\not=\emptyset$ and $\sigma_R=\emptyset$, then
\begin{enumerate}
\item[---]  if $(\hat\sigma)>0$, then $\mathsf{smax}(\sigma)=\hat{\sigma}$.
\item[---]   if $(\hat\sigma)<0$, then $\mathsf{smax}(\sigma)=\mathsf{smax}(\sigma_L)$.
\end{enumerate}

\item If $\sigma_L\not=\emptyset$ and $\sigma_R\not=\emptyset$ and $\hat{\sigma}>0$, then
\begin{enumerate}
\item[---]   if $\min(|\sigma_L|)>\min(|\sigma_R|)$, then $\mathsf{smax}(\sigma)=\mathsf{smax}(\sigma_L)$.
\item[---]   if $\min(|\sigma_L|)<\min(|\sigma_R|)$, then $\mathsf{smax}(\sigma)=\mathsf{smax}(\sigma_R)$.
\end{enumerate}

\item If $\sigma_L\not=\emptyset$ and $\sigma_R\not=\emptyset$ and $\hat{\sigma}<0$, then
\begin{enumerate}
\item[---]  if $\min(|\sigma_L|)<\min(|\sigma_R|)$, then $\mathsf{smax}(\sigma)=\mathsf{smax}(\sigma_L)$.
\item[---]  if $\min(|\sigma_L|)>\min(|\sigma_R|)$, then $\mathsf{smax}(\sigma)=\mathsf{smax}(\sigma_R)$.
\end{enumerate}

\end{enumerate}

Note that for the flip equivalence class $[\sigma]$, $\mathsf{smax}([\sigma]):=\mathsf{smax}(\sigma)$ is well-defined as $\mathsf{smax}(\sigma')=\mathsf{smax}(\sigma)$ for all $\sigma'\in[\sigma]$.
For example, we have
$\mathsf{smax}(27\bar{8}16\bar{9}\bar{3}\bar{4}5)=5$ 
as the above process to determine the signed maximum
results in $27\bar{8}16\bar{9}\bar{3}\bar{4}5\to 6\bar{9}\bar{3}\bar{4}5 \to \bar{4}5 \to 5$. Another example
is $\mathsf{smax}([2\bar{1}3\bar{4}])=2$ as one can check
that all elements in $[2\bar{1}3\bar{4}]=\{2\bar{1}3\bar{4}, 
2\bar{1}\bar{4}3, 3\bar{4}\bar{1}2, \bar{4}3\bar{1}2\}$ have the same signed maximum.

\medskip

Let $\mathcal{FL}^{(B)}_n$ be the set of flip equivalence classes of signed permutations in $\mathfrak{B}_n$ whose signed maximums are positive. 
Also define $$\mathcal{FL}^{(B)}_{n,k}:=\{[\sigma]\in\mathcal{FL}^{(B)}_n,\mathsf{smax}([\sigma])=k\}$$
for $1\leq k\leq n$. We have
\begin{align*}
    \mathcal{FL}^{(B)}_1=&\{[1]\};\\
   \mathcal{FL}^{(B)}_2=&\{[12],[1\bar{2}],[\bar{1}2]\};\\
    \mathcal{FL}^{(B)}_3=&\{[123],[213],[\bar{1}23],[1\bar{2}3],[12\bar{3}],[\bar{1}\bar{2}3],[\bar{1}2\bar{3}],[1\bar{2}\bar{3}],[2\bar{1}3],[\bar{2}13],[2\bar{1}\bar{3}]\}.
  \end{align*}

Similarly, let $\mathcal{FL}^{(D)}_n$ be the set of flip equivalence classes of signed permutations in $\mathfrak{B}_n$ whose signed maximums are negative, 
and define $$\mathcal{FL}^{(D)}_{n,k}:=\{[\sigma]\in\mathcal{FL}^{(D)}_n,\mathsf{smax}([\sigma])=-k\}$$ for $1\leq k\leq n$.
We have
\begin{align*}
      \mathcal{FL}^{(D)}_1=&\{[\bar{1}]\};\\
   \mathcal{FL}^{(D)}_2=&\{[\bar{1}\bar{2}]\};\\
    \mathcal{FL}^{(D)}_3=&\{[\bar{1}\bar{2}\bar{3}],[\bar{2}1\bar{3}],[\bar{2}\bar{1}3],[2\bar{1}\bar{3}],[\bar{2}\bar{1}\bar{3}]\}.
\end{align*}

We will see that both are counted by Springer numbers of types $B_n$ and $D_n$. Moreover, define the statistic {\it signed peak} $\mathsf{spk}$ for $\sigma\in \mathfrak{B}_n$ by $$\mathsf{spk}(\sigma):=\#\{i:1\leq i\leq n,\sigma_i<0,|\sigma_{i-1}|<|\sigma_i|>|\sigma_{i+1}|\}$$ with the convention $\sigma_0=\sigma_{n+1}=0$. Note that $\mathsf{spk}([\sigma]):=\mathsf{spk}(\sigma)$ is well-defined as $\mathsf{spk}(\sigma')=\mathsf{spk}(\sigma)$ for all $\sigma'\in[\sigma]$. 
With this we obtain the third realization of the refined Arnold family. 

\begin{thm}\label{flpoly} (\textbf{Main Theorem III})
For $1\leq k\leq n$, we have
$$V_{n,k}(t)=\sum_{\sigma\in\mathcal{FL}^{(B)}_{n,n-k+1}}t^{n+1-2\cdot\mathsf{spk}(\sigma)}$$
and $$V_{n,-k}(t)=\sum_{\sigma\in\mathcal{FL}^{(D)}_{n,n-k+1}}t^{n+1-2\cdot \mathsf{spk}(\sigma)}.$$
\end{thm}
In other words, $$(\mathcal{FL}^{(D)}_{n,1},\ldots,\mathcal{FL}^{(D)}_{n,n},\mathcal{FL}^{(B)}_{n,n},\ldots,\mathcal{FL}^{(B)}_{n,1})$$ 
together with the statistic $\mathsf{spk}$ forms a refined Arnold family.

\subsection{complete increasing binary trees with $n$ labelled nodes}

\begin{figure}[ht]
\begin{center}
\includegraphics[width=5.4in]{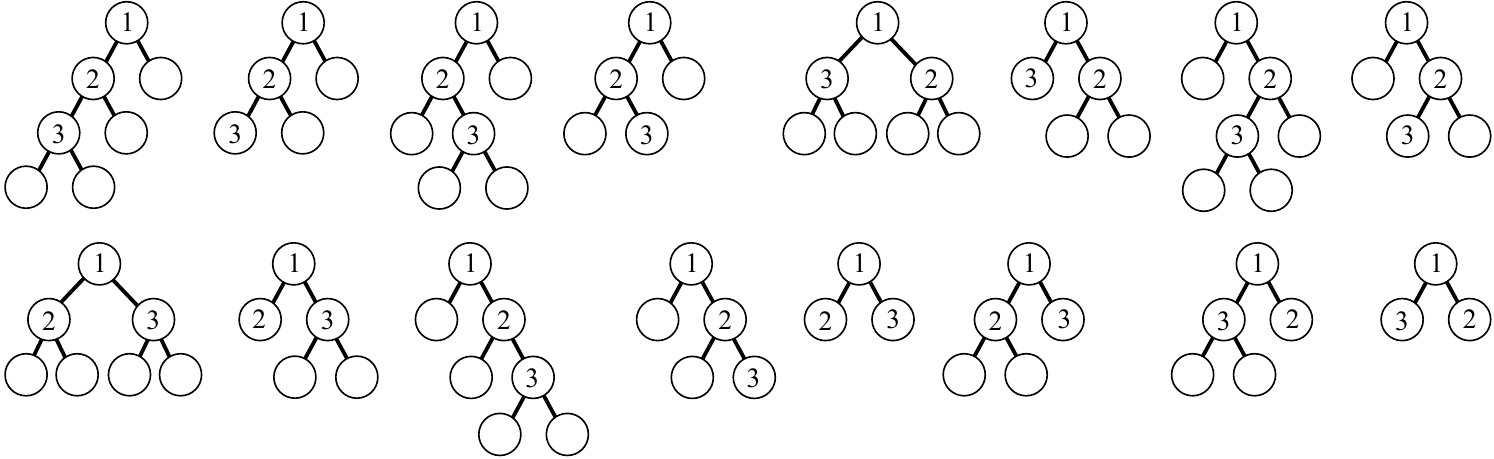}
\end{center}
\caption{\small The complete increasing binary trees with three labelled nodes.}
\label{fig:plane-labelled-binary-trees}
\end{figure}

A (plane) binary tree is \emph{complete} if every non-leaf node has two children. 
Let $\T_n$ be the set of complete binary trees such that there are $n$ nodes labelled by integers from $1$ to $n$ but some leaves can be left unlabelled (called \emph{empty leaves}), and along the path from the root to any leaf the labels are increasing. These trees are first introduced by~\cite{JosuatVerges_10}, and in the rest of the paper, by \emph{complete increasing binary trees} we always mean the trees from the set $\T_n$. Figure \ref{fig:plane-labelled-binary-trees} shows the sixteen complete increasing binary trees with three labelled nodes. 

Recently a realization of the refined Arnold family in terms of complete increasing binary trees is found by Eu and Fu~\cite{Eu_Fu_23}.
Given $\tau\in\mathcal{T}_n$, the {\it rightmost path} of $\tau$ is a sequence of nodes $(v_1,v_2,\ldots,v_d)$ where $v_1$ is the root and $v_{i+1}$ is the right child of $v_i$ for all $i=1,2,\ldots,d-1$. The node $v_d$ is called the rightmost leaf of $\tau$. Let $\mathcal{T}^\circ_n$ ($\mathcal{T}^*_n$, respectively) denote the subset of $\mathcal{T}_n$ consists of the trees whose rightmost leaf is empty (labelled, respectively). The label of $v_d$ ($v_{d-1}$, respectively) is called the {\it rightmost label} of $\tau$ for $\tau\in\mathcal{T}^\circ_n$ ($\tau\in\mathcal{T}^*_n$, respectively). For $1\leq k\leq n$, let $\mathcal{T}^\circ_{n,k}\subset\mathcal{T}^\circ_n$ ($\mathcal{T}^*_{n,k}\subset\mathcal{T}^*_n$, respectively) be the subset of trees whose rightmost label is $k$.
 Let $\mathsf{emp}(\tau)$ be the number of empty leaves in $\tau$, we have the following.
\begin{thm}[\cite{Eu_Fu_23}]\label{treeV}
For $1\leq k\leq n$, $$V_{n,k}(t)=\sum_{\tau\in\mathcal{T}^\circ_{n,n-k+1}}t^{\mathsf{emp}(\tau)}$$ and $$V_{n,-k}(t)=\sum_{\tau\in\mathcal{T}^*_{n,n-k+1}}t^{\mathsf{emp}(\tau)}.$$
\end{thm}
It is then natural to seek bijection between complete increasing binary trees to cycle-up-down permutations, valley signed permutations, and flip equivalence classes.

\bigskip

The main results of this paper are twofold. Firstly, we will give definitions (if needed) of types $B_n$ and $D_n$ versions of the above structures and find suitable statistics such that each is a new realization of the refined Arnold family. Secondly, as the set of complete increasing binary trees is the key structure for linking other Arnold families together (see~\cite{Eu_Fu_23}), we will also give bijections between these new combinatorial objects with the complete increasing binary trees. 

The rest of the paper is organized as follows. In Section 3 we give a realization of the refined Arnold family in terms of cycle-up-down permutations and in Section 4 in terms of valley signed permutations. 
The interlude Sections 5 (resp. Section 6) are devoted to bijections 
between cycle-up-down permutations (resp. valley signed permutations) with complete increasing binary trees. In Section 7 we investigate Knuth's flip equivalence classes on signed permutations.

\section{Proofs of Main Theorem I}
For $1\leq k\leq n$, define the polynomials 
$$CUD^{(B)}_{n,k}=CUD^{(B)}_{n,k}(t):=\sum_{\sigma\in\mathcal{CUD}^{(B)}_{n,k}}t^{n+1-2\cdot \mathsf{npk}(\sigma)}$$ 
and $$CUD^{(D)}_{n,k}=CUD^{(D)}_{n,k}(t):=\sum_{\sigma\in\mathcal{CUD}^{(D)}_{n,k}}t^{n+1-2\cdot \mathsf{npk}(\sigma)}.$$
To prove Theorem~\ref{cudpoly}, by referring to the recurrences of $V_{n,k}(t)$
we need to prove the following three recurrences:
\begin{enumerate}
\item $CUD^{(D)}_{n,k}=CUD^{(D)}_{n,k-1}+t^{-1}CUD^{(B)}_{n-1,k-1}$  for $1<k\leq n$;
\item $CUD^{(B)}_{n,n}=t^2CUD^{(D)}_{n,n}$ for $n\ge 2$;
\item $CUD^{(B)}_{n,k}=CUD^{(B)}_{n,k+1}+tCUD^{(D)}_{n-1,k}$ for $1\leq k<n$.
\end{enumerate}
We investigate them in the propositions below.  
\begin{prop}\label{cudrec1}
For $1<k\leq n$, we have 
$$CUD^{(D)}_{n,k}=CUD^{(D)}_{n,k-1}+t^{-1}CUD^{(B)}_{n-1,k-1}.$$ 
\end{prop}
\begin{proof}
We shall construct a bijection 
$$\psi: \mathcal{CUD}^{(D)}_{n,k} \longrightarrow \mathcal{CUD}^{(D)}_{n,k-1}\cup\mathcal{CUD}^{(B)}_{n-1,k-1}.$$ 

Given $\sigma=(a_{1,1},a_{1,2},\ldots)\ldots(a_{m-1,1},a_{m-1,2},\ldots)(k,\bar{k})\in\mathcal{CUD}^{(D)}_{n,k}$.

(i) If $a_{m-1,1}=k-1$.  We first delete the cycle $(k,\bar{k})$ and let 
$$\psi(\sigma)=(a_{1,1}',a_{1,2}',\ldots)(a_{2,1}',a_{2,2}',\ldots)\ldots(a_{m-1,1}',a_{m-1,2}',\ldots),$$ 
where $a_{i,j}'=\begin{cases}
    a_{i,j}, &\mbox{if }|a_{i,j}|<k,\\
    a_{i,j}-1, &\mbox{if }a_{i,j}>k,\\
    a_{i,j}+1, &\mbox{if }a_{i,j}<-k.
    \end{cases}$.
It is clear that $\psi(\sigma)\in\mathcal{CUD}^{(B)}_{n-1,k-1}$.
    
    As for the polynomial, since $\mathsf{npk}(\psi(\sigma))=\mathsf{npk}(\sigma)-1$, so 
$$(n-1)+1-2\mathsf{npk}(\psi(\sigma))=(n+1-2\mathsf{npk}(\sigma))+1$$ and we get the term $t^{-1}CUD^{(B)}_{n-1,k-1}$.

(ii) If $a_{m-1,1}\not=k-1$. We let 
$$\psi(\sigma)=(a_{1,1}',a_{1,2}',\ldots)(a_{2,1}',a_{2,2}',\ldots)\ldots(a_{m-1,1}',a_{m-1,2}',\ldots),(k-1, \overline{k-1})$$
 where $a_{i,j}'=\begin{cases}
    a_{i,j}, &\mbox{if }|a_{i,j}|\mbox{ is not }k\mbox{ or }k-1,\\
    \mathsf{sgn}(a_{i,j})k, &\mbox{if }|a_{i,j}|=k-1,\\
    \end{cases}$

It is also clear that $\psi(\sigma)\in\mathcal{CUD}^{(D)}_{n,k-1}$ and $\mathsf{npk}(\psi(\sigma))=\mathsf{npk}(\sigma)$. 
Hence we have the term $CUD^{(D)}_{n,k-1}$.

The inverse map $\psi^{-1}$ can be established directly by the reversing the above process the lemma is proved.
\end{proof}

For examples, $\sigma=(1,\bar{5},\bar{2})(3,4)(6,9,\bar{8})(7,\bar{7})\in\mathcal{CUD}^{(D)}_{9,7}$ is mapped to $$\psi(\sigma)=(1,\bar{5},\bar{2})(3,4)(6,8,\bar{7})\in \mathcal{CUD}^{(B)}_{8,6}$$ and $\sigma'=(1,\bar{5},\bar{2})(3,\bar{6})(4,9,\bar{8})(7,\bar{7})\in\mathcal{CUD}^{(D)}_{9,7}$ mapped to
$$\psi(\sigma')=(1,\bar{5},\bar{2})(3,\bar{7})(4,9,\bar{8})(6,\bar{6})\in \mathcal{CUD}^{(D)}_{9,6}.$$

\begin{prop}\label{cudrec2}
 For $n\geq 2$, we have $CUD^{(B)}_{n,n}=t^2CUD^{(D)}_{n,n}$.
\end{prop}
\begin{proof}
The bijection $\psi: \mathcal{CUD}^{(B)}_{n,n} \longrightarrow \mathcal{CUD}^{(D)}_{n,n}$ is easy:
we replace the last cycle $(n)$ in $\sigma \in  \mathcal{CUD}^{(B)}_{n,n}$ by the cycle $(n,\bar{n})$ (and vice versa). 
Notice that $\mathsf{npk}(\psi(\sigma))=\mathsf{npk}(\sigma)+1$ and 
$n+1-2\mathsf{npk}(\psi(\sigma))=(n+1-2\mathsf{npk}(\sigma))-2$, we get the term 
$t^2CUD^{(D)}_{n,n}$ as desired.
\end{proof}


\begin{prop}\label{cudrec3}
For $1\leq k<n$, we have $$CUD^{(B)}_{n,k}=CUD^{(B)}_{n,k+1}+tCUD^{(D)}_{n-1,k}.$$
\end{prop}

\begin{proof}
We construct a bijection $$\psi: \mathcal{CUD}^{(B)}_{n,k} \longrightarrow \mathcal{CUD}^{(B)}_{n,k+1}\cup\mathcal{CUD}^{(D)}_{n-1,k}.$$ 
Given $\sigma \in\mathcal{CUD}^{(B)}_{n,k}$. Suppose the first entry of the last cycle of is $k$, i.e., it is of the form
$\sigma=(a_{1,1},a_{1,2},\ldots)\ldots(a_{m-1,1},a_{m-1,2},\ldots)(k,\ldots)$. We have three cases.

(i) If the last cycle of $\sigma$ is $(k,\overline{k+1})$, let $$\psi(\sigma)=(a_{1,1}',a_{1,2}',\ldots)(a_{2,1}',a_{2,2}',\ldots)\ldots(k,\bar{k}),$$ where

 $$a_{i,j}'=\begin{cases}
    a_{i,j}, &\mbox{if }|a_{i,j}|<k,\\
    a_{i,j}-1, &\mbox{if }a_{i,j}>k+1,\\
    a_{i,j}+1, &\mbox{if }a_{i,j}<-(k+1).
    \end{cases}$$

It is clear that $\psi(\sigma)\in\mathcal{CUD}^{(D)}_{n-1,k}$ and $\mathsf{npk}(\psi(\sigma))=\mathsf{npk}(\sigma)$. 
Since $(n-1)+1-2\mathsf{npk}(\psi(\sigma))=(n+1-2\mathsf{npk}(\sigma))-1$, we get the term $tCUD^{(D)}_{n-1,k}$.

(ii) If the last cycle of $\sigma$ is not $(k,\overline{k+1})$ but there exist $|a_{m,\ell}|=k+1$ for some $\ell\geq 2$, then we divide the last cycle into two cycles. That is, let $$\psi(\sigma)=(a_{1,1},a_{1,2},\ldots)\ldots(k,a_{m,2},\ldots , a_{m,\ell-1})(k+1,a_{m,\ell+1},\ldots)$$ and we have $\psi(\sigma)\in\mathcal{CUD}^{(B)}_{n,k+1}$ and $\mathsf{npk}(\psi(\sigma))=\mathsf{npk}(\sigma)$.

(iii) If the last cycle doesn't contain $k+1$ or $\overline{k+1}$, let 
$$\psi(\sigma)=(a_{1,1}',a_{1,2}',\ldots)(a_{2,1}',a_{2,2}',\ldots)\ldots(a_{m,1}',a_{m,2}',\ldots),$$ 
where $a_{i,j}'=\begin{cases}
    a_{i,j}, &\mbox{if }|a_{i,j}|\mbox{ is neither }k\mbox{ nor }k-1,\\
    \mathsf{sgn}(a_{i,j})k, &\mbox{if }|a_{i,j}|=k-1,\\
    \mathsf{sgn}(a_{i,j})(k-1), &\mbox{if }|a_{i,j}|=k.
    \end{cases}$\\ 
Then $\psi(\sigma)\in\mathcal{CUD}^{(B)}_{n,k+1}$ and $\mathsf{npk}(\psi(\sigma))=\mathsf{npk}(\sigma)$.

Note that cases (ii) and (iii) lead to the term $CUD^{(B)}_{n,k+1}$. The inverse map $\psi^{-1}$ can be established by reversing the process, 
hence the recurrence is verified.
\end{proof}

For examples, $\sigma=(1,\bar{3},\bar{2})(4)(5,\bar{6})(7,9,\bar{8})\in\mathcal{CUD}^{(B)}_{9,7}$ is
mapped to $$\psi(\sigma)=(1,\bar{3},\bar{2})(4)(5,\bar{6})(7,9)(8)\in \mathcal{CUD}^{(B)}_{9,8}$$  by case (ii),
while by case (iii), $\sigma'=(1,\bar{3},\bar{2})(4)(5,\bar{8},\bar{6})(7,9)\in\mathcal{CUD}^{(B)}_{9,7}$ is mapped to 
$$\psi(\sigma')=(1,\bar{3},\bar{2})(4)(5,\bar{7},\bar{6})(8,9)\in \mathcal{CUD}^{(B)}_{9,8}.$$

We are ready to prove Main Result I.
\begin{proof}[Proof of Theorem~\ref{cudpoly}]
It remains to check the initial cases. For $n=1$, we have $\mathsf{npk}((1))=0$ and $\mathsf{npk}((1,\bar{1}))=1$. 
Hence $$CUD^{(B)}_{1,1}(t)=t^{1+1-2\mathsf{npk}((1))}=t^2=V_{1,1}(t)$$
and
$$CUD^{(D)}_{1,1}(t)=t^{1+1-2\mathsf{npk}((1,\bar{1}))}=1=V_{1,-1}(t).$$
Since $\mathcal{CUD}^{(D)}_{n,1}$ are empty for all $n\geq 2$, we have
$$CUD^{(D)}_{n,1}(t)=0=V_{n,-n}(t)\quad(n\geq 2).$$
Together with Propositions~\ref{cudrec1}, \ref{cudrec2}, and \ref{cudrec3}, the proof is completed.
\end{proof}

\section{Proof of Main Theorem II}
In this section we prove Theorem 2.2. 
For $1\leq k\leq n$, define 
$$VS^{(B)}_{n,k}=VS^{(B)}_{n,k}(t):=\sum_{\sigma\in\mathcal{VS}^{(B)}_{n,k}}t^{n+1-2\mathsf{neg}(\sigma)}$$
and $$VS^{(D)}_{n,k}=VS^{(D)}_{n,k}(t):=\sum_{\sigma\in\mathcal{VS}^{(D)}_{n,k}}t^{n+1-2\mathsf{neg}(\sigma)}.$$
We are going to show that both have the same recurrence relations as $V_{n,k}$.

\begin{prop}\label{vsrec}
The following recurrence relations hold.
\begin{enumerate}
    \item $VS^{(D)}_{n,k}=VS^{(D)}_{n,k-1}+t^{-1}VS^{(B)}_{n-1,k-1}$, for $1<k\leq n$. 
    \item $VS^{(B)}_{n,n}=t^2VS^{(D)}_{n,n}$, for $n\geq 2$. 
    \item $VS^{(B)}_{n,k}=VS^{(B)}_{n,k+1}+tVS^{(D)}_{n-1,k}$, for $1\leq k<n$. 
\end{enumerate} 
\end{prop}
\begin{proof}
(i) We  construct a bijection $$\psi: \mathcal{VS}^{(D)}_{n,k}\longrightarrow \mathcal{VS}^{(D)}_{n,k-1}\cup\mathcal{VS}^{(B)}_{n-1,k-1}.$$ 
Given $\sigma=\sigma_1\sigma_2\ldots \sigma_n\in\mathcal{VS}^{(D)}_{n,k}$ with $\sigma_1=-k$.

Case 1: If $\sigma_2=k-1$. We first remove $\sigma_1$ and let $\psi(\sigma)=\sigma_1'\sigma_2'\ldots \sigma_{n-1}'$, where 
$$\sigma_i'=\begin{cases}
    \sigma_{i+1}, &\mbox{if }|\sigma_{i+1}|<k,\\
    \sigma_{i+1}-1, &\mbox{if }\sigma_{i+1}>k,\\
    \sigma_{i+1}+1, &\mbox{if }\sigma_{i+1}<-k.
    \end{cases}$$
It is clear that $\psi(\sigma)\in\mathcal{VS}^{(B)}_{n-1,k-1}$ and $\mathsf{neg}(\psi(\sigma))=\mathsf{neg}(\sigma)-1$. 
Since $(n-1)+1-2\mathsf{neg}(\psi(\sigma))=(n+1-2\mathsf{neg}(\sigma))+1$, we get the term $t^{-1}VS^{(B)}_{n-1,k-1}$.

Case 2: If $\sigma_2\not=k-1$, let $\psi(\sigma)=\sigma_1'\sigma_2'\ldots \sigma_n'$ with
$$\sigma_i'=\begin{cases}
    \sigma_i, &\mbox{if }|\sigma_i|\mbox{ is neither }k\mbox{ nor }k-1,\\
    \mathsf{sgn}(\sigma_i)k, &\mbox{if }|\sigma_i|=k-1,\\
    \mathsf{sgn}(\sigma_i)(k-1), &\mbox{if }|\sigma_i|=k.
    \end{cases}$$ 
Then we have $\psi(\sigma)\in\mathcal{VS}^{(D)}_{n,k-1}$, $\mathsf{neg}(\psi(\sigma))=\mathsf{neg}(\sigma)$ and we obtain the term $VS^{(D)}_{n,k-1}$.

(ii) There is an immediate bijection $\psi: \mathcal{VS}^{(B)}_{n,n}\longrightarrow \mathcal{VS}^{(B)}_{n,n}$.
For $\sigma\in \mathcal{VS}^{(B)}_{n,n}$, the image $\psi(\sigma)$ is obtained by replacing $\sigma_1(=n)$ by $\bar{n}$.
Since $\mathsf{neg}(\psi(\sigma))=\mathsf{neg}(\sigma)+1$, we have $n+1-2\mathsf{neg}(\psi(\sigma))=(n+1-2\mathsf{neg}(\sigma))-2$ and recurrence (ii) is proved.

(iii) We construct a bijection $$\psi: \mathcal{VS}^{(B)}_{n,k} \longrightarrow \mathcal{VS}^{(B)}_{n,k+1}\cup\mathcal{VS}^{(D)}_{n-1,k}.$$ 

Given $\sigma\in\mathcal{VS}^{(B)}_{n,k}$ with $\sigma_1=k$.

Case 1. If $\sigma_2=\overline{k+1}$. We firstly remove $\sigma_1$ and let $\psi(\sigma)=\sigma_1'\sigma_2'\ldots \sigma_{n-1}'$, where 
$$\sigma_i'=\begin{cases}
    \sigma_{i+1}, &\mbox{if }|\sigma_{i+1}|<k,\\
    \sigma_{i+1}-1, &\mbox{if }\sigma_{i+1}>k.\\
    \sigma_{i+1}+1, &\mbox{if }\sigma_{i+1}<-k.
    \end{cases}$$
It can be seen that $\psi(\sigma)\in\mathcal{VS}^{(D)}_{n-1,k}$ and $\mathsf{neg}(\psi(\sigma))=\mathsf{neg}(\sigma)$. 
Since $(n-1)+1-2\mathsf{neg}(\psi(\sigma))=(n+1-2\mathsf{neg}(\sigma))-1$, we get the term $tVS^{(D)}_{n-1,k}$. 

Case 2. If $\sigma_2\not=\overline{k+1}$, then we let $\psi(\sigma)=\sigma_1'\sigma_2'\ldots\sigma_n'$ with
$$\sigma_i'=\begin{cases}
    \sigma_i, &\mbox{if }|\sigma_i|\mbox{ is neither }k\mbox{ nor }k+1,\\
    \mathsf{sgn}(\sigma_i)k, &\mbox{if }|\sigma_i|=k+1,\\
    \mathsf{sgn}(\sigma_i)(k+1), &\mbox{if }|\sigma_i|=k.
    \end{cases}$$ 
Now $\psi(\sigma)\in\mathcal{VS}^{(B)}_{n,k+1}$ and $\mathsf{neg}(\psi(\sigma))=\mathsf{neg}(\sigma)$, which gives us the term $VS^{(B)}_{n,k+1}$. 
The inverse map $\psi^{-1}$ can be established by reversing the operations and (iii) is verified.
\end{proof}

For examples, by case 2 of (i), 
$\sigma=\bar{7}4281\bar{5}3\bar{9}(10)6\in\mathcal{VS}^{(D)}_{10,7}$ is mapped to 
$\bar{6}4281\bar{5}3\bar{9}(10)7 \in \mathcal{VS}^{(D)}_{10,6}$.
Also, by the case 2 of (iii), $\sigma=7985\bar{6}41\bar{3}2\in\mathcal{VS}^{(B)}_{9,7}$
is mapped to  $8975\bar{6}41\bar{3}2 \in \mathcal{VS}^{(B)}_{9,8}$.\\

Now we are ready to prove Main Result II.
\begin{proof}[Proof of Theorem~\ref{vspoly}]
It remains to check the initial values. For the case $n=1$, we have $\mathsf{neg}(1)=0$ and $\mathsf{neg}(\bar{1})=1$. Hence $$VS^{(B)}_{1,1}(t)=t^{1+1-2\mathsf{neg}(1)}=t^2=V_{1,1}(t)$$
and
$$VS^{(D)}_{1,1}(t)=t^{1+1-2\mathsf{neg}(\bar{1})}=1=V_{1,-1}(t).$$
Since $\mathcal{VS}^{(D)}_{n,1}$ are empty for all $n\geq 2$, we have
$$VS^{(D)}_{n,1}(t)=0=V_{n,-n}(t)$$
for $n\geq 2$. Together with Proposition~\ref{vsrec}, the proof is completed.
\end{proof}

\section{Complete increasing binary trees and Cycle-up-down permutations}
Recall that $\T_n$ is the set of complete increasing binary trees with $n$ labelled nodes 
and $\mathcal{T}^\circ_n$ ($\mathcal{T}^*_n$, respectively) are those trees with empty (resp. labelled) rightmost leaf.
In this section, we shall establish the bijections 
$$\phi^{(B)}_C:\mathcal{CUD}^{(B)}_n\to\mathcal{T}^\circ_n \quad \text{and} \quad \phi^{(D)}_C:\mathcal{CUD}^{(D)}_n\to\mathcal{T}^*_n.$$
The subscript $C$ indicates that it is for cycle-up-down permutations.

\subsection{Type $B_n$}
We roughly describe the flowchart of the desired bijection. We map an up-down sequence to a non-plane complete increasing tree by
{\bf Algorithm 1}. By {\bf Algorithm 2} we map each cycle of $\sigma\in\mathcal{CUD}^{(B)}_n$ to a complete increasing binary tree with the help of {\bf Algorithm 1}. By combing corresponding trees from all cycles we reach our desired bijection. Our running example is $\sigma=(1,\bar{3},\bar{2})(4)(5,\bar{6})(7,9,\bar{8})\in\mathcal{CUD}^{(B)}_{9,7}$. 
We need two operations first.
\begin{enumerate}
    \item {\it Double bracket decomposition:} Decompose a sequence $a_1,a_2,\ldots,a_n$ of distinct integer into $((a_1,\ldots,a_{i-1})a_i(a_{i+1},\ldots,a_n))$, where $a_i$ is the smallest entry.
    \item {\it Complement:} Replace the $i$-th smallest entry of a sequence $\mathcal{A}=a_1,a_2,\ldots,a_n$ of distinct integers 
by the $i$-th greatest entry for each $1\le i\le n$. Denote the resulting sequence by $\mathcal{A}^c$. 
\end{enumerate}

{\bf Algorithm 1} is actually adapted from~\cite{Donaghey_75}, which maps a sequence of distinct integers into a non-plane complete increasing binary tree.\\

\noindent \textbf{Algorithm 1:}\\
{\bf Input :} A sequence $\mathcal{A}=a_1,a_2,\ldots,a_n$ of distinct integers.

Set $\mathcal{T}(\mathcal{A})$ to be an empty tree. Apply the following procedures recursively.
\begin{enumerate}
    \item If $\max(\mathcal{A})$ is on the left of $\min(\mathcal{A})$ : Let $\mathcal{A}:=\mathcal{A}^c$ and go to {(ii)}
    \item If $\max(\mathcal{A})$ is on the right of $\min(\mathcal{A})$ : Apply double bracketing decomposition on $\mathcal{A}$ and write $\mathcal{A}=((\mathcal{A}_L)a_i(\mathcal{A}_R))$. Construct a tree rooted at $a_i$ and define the children of $a_i$ to be two nodes labelled by $\min(\mathcal{A}_L)$ and $\min(\mathcal{A}_R)$ (If any of the sequence is empty, the child is defined to be empty). 
	\item Apply the above procedure recursively on $\mathcal{A}_L$ and $\mathcal{A}_R$.
\end{enumerate}
 {\bf Output :} 
A complete increasing binary (non-plane) tree $\mathcal{T}(\mathcal{A})$. Note that the leaves of $\mathcal{T}(\mathcal{A})$ are all empty.\\

\begin{figure} [htb]
\begin{center}
\includegraphics[scale=0.2]{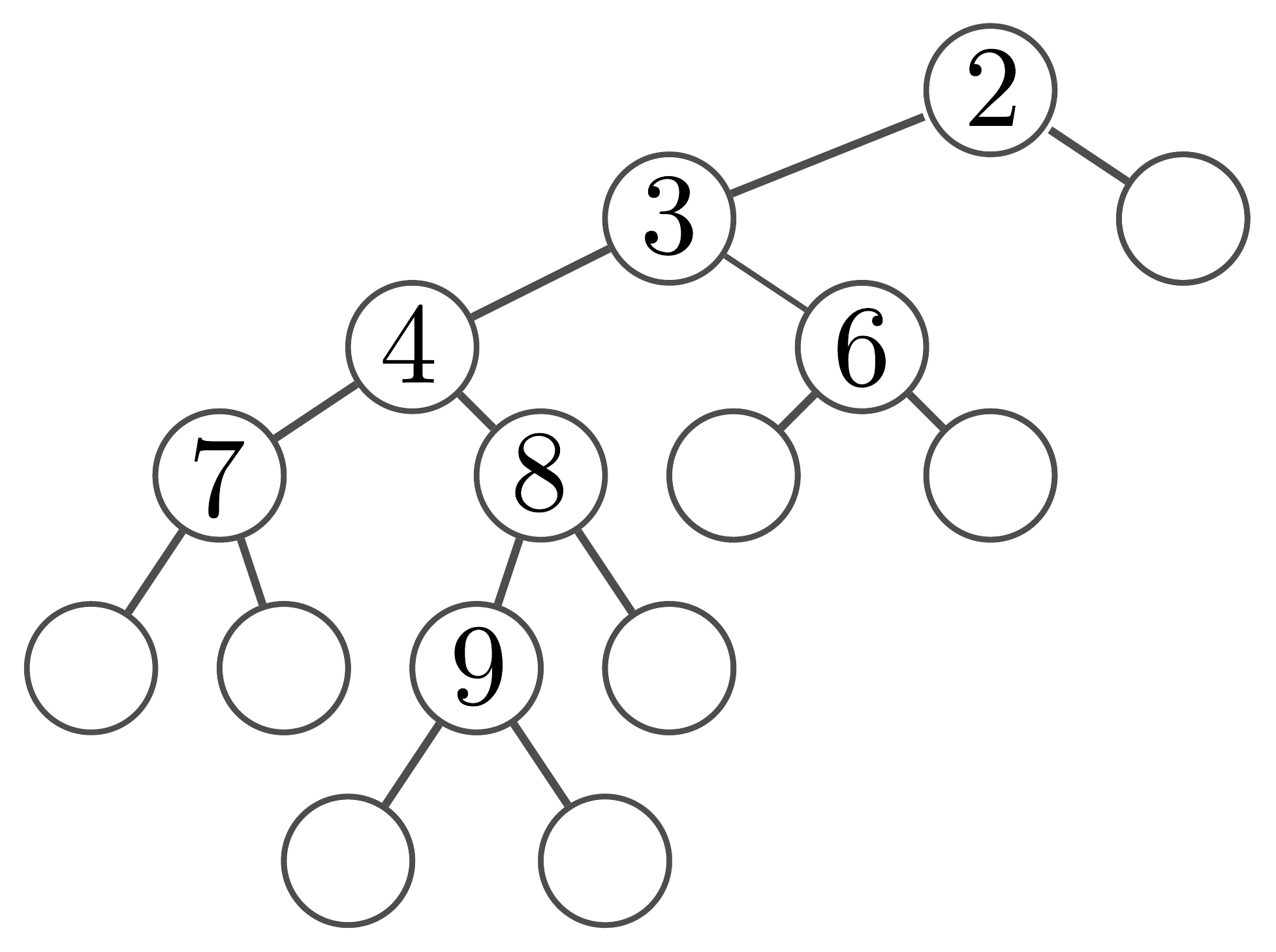}
\end{center}
\caption{\small From a sequence $\mathcal{A}$ to the complete increasing binary tree $T(\mathcal{A})$.}
 \label{algo1}
\end{figure}

For example, Figure \ref{algo1} maps $\mathcal{A}=2,6,3,7,4,9,8$ to the 
corresponding tree $T(\mathcal{A})$ with these numbers as the labels.

We say a permutation $\sigma=\sigma_1,\sigma_2,\ldots,\sigma_{|A|}$ of elements of a set $A$ of distinct integers {\it up-down} if
$\sigma_1<\sigma_2>\sigma_3<\cdots$. Let $\mathcal{UD}_A$ be the set of up-down permutations
and $\mathcal{UD}_{A,1}\subseteq \mathcal{UD}_A$ consisting of those with $\sigma_1=\min(\sigma)$.
The following is an easy consequence of {\bf Algorithm 1}.

\begin{cor} {\bf Algorithm 1} induces a bijection between 
\begin{enumerate} 
\item  $\mathcal{UD}_A$ and the set of complete increasing (non-plane) tree with nodes labelled by elements of $A$ such that all leaves are empty.
\item  $\mathcal{UD}_{A,1}$ and the set of complete increasing (non-plane) tree with nodes labelled by elements of $A$ such that all leaves are empty and the root has an empty child. 
\end{enumerate}
\end{cor}

Given $\sigma=(a_{1,1},a_{1,2},\ldots)(a_{2,1},a_{2,2},\ldots)\ldots(a_{m,1},a_{m,2},\ldots)\in\mathcal{CUD}^{(B)}_n$. 
For the cycle $C_i=(a_{i,1},a_{i,2},\ldots)$, $1\le i\le m$, we let $\mathsf{int}(C_i)=\{|a_{i,1}|,|a_{i,2}|,\ldots\}$. 
On the other hand, for a finite set $A$ of positive integers we let $\mathcal{T}^{(1)}_A$ denote the set of complete increasing binary plane trees with empty leaves, satisfying that (i) the nodes are labelled using the integers in $A$, and (ii) the rightmost labelled node is the root. 
The goal of {\bf Algorithm 2} is to map a cycle $C_i$ to a tree $\tau\in\mathcal{T}^{(1)}_{\mathsf{int}(C_i)}$.\\ 

\noindent \textbf{Algorithm 2}:\\
\noindent \textbf{Input :} A cycle $C_i=(a_{i,1},a_{i,2},\ldots,a_{i,d})$ of a signed cycle-up-down permutation $\sigma$.

First we let $\mathcal{A}=|a_{i,1}|, |a_{i,2}|,\ldots , |a_{i,d}|$ and construct $\mathcal{T}(\mathcal{A})$ using {\bf Algorithm 1}.

For $1\leq j\leq d$, do the followings.
\begin{enumerate}
\item  If the node $v$ with label $|a_{i,j}|$ contains two empty children in $\mathcal{T}(\mathcal{A})$: 
	        \begin{itemize}
            \item[---] If $a_{i,j}<0$, then remove the two empty children.
        \end{itemize}
\item  If the node $v$ with label $|a_{i,j}|$ contains exactly one labelled child in $\mathcal{T}(\mathcal{A})$:
        \begin{itemize}
            \item[---] If $a_{i,j}<0$, let the labelled child be the right child of $v$ and the empty node be the left child of $v$.  
            \item[---] If $a_{i,j}>0$, let the labelled child be the left child of $v$ and the empty node be the right child of $v$.
        \end{itemize}
\item If the node $v$ with label $|a_{i,j}|$ contains two labelled children in $\mathcal{T}(\mathcal{A})$:
        \begin{itemize}
            \item[---] If $a_{i,j}<0$, let the child with smaller label be the right child of $v$ and the child with larger label be the left child of $v$.  
            \item[---] If $a_{i,j}>0$, let the child with smaller label be the left child of $v$ and the child with larger label be the right child of $v$. 
        \end{itemize}
\end{enumerate}
\noindent \textbf{Output :} The resulting tree $\tau_i$. 
Note that it is easy to see that $\tau_i \in \mathcal{T}^{(1)}_{\mathsf{int}(C_i)}$.

\medskip

\begin{figure}[htb]
\begin{center}
\includegraphics[scale=0.2]{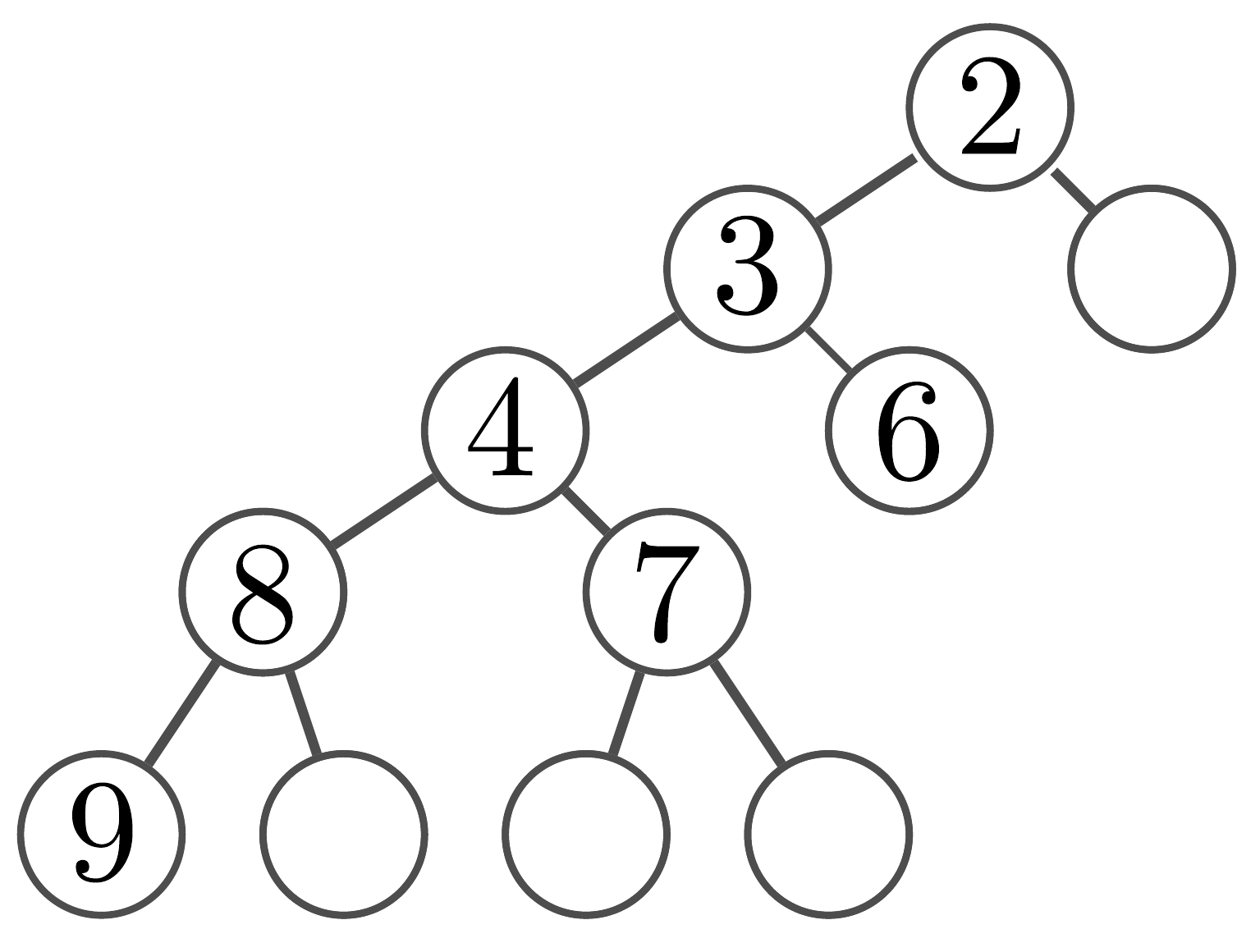}
\end{center}
\caption{\small The tree in $\mathcal{T}^{(1)}_{\mathsf{int}(C)}$ from the cycle
$C=(2,\bar{6},3,7,\bar{4},\bar{9},8)$.}
\label{algo2}
\end{figure}

For example, the cycle $C=(2,\bar{6},3,7,\bar{4},\bar{9},8)$
maps to a tree in $\mathcal{T}^{(1)}_{\mathsf{int}(C)}$ as in Figure~\ref{algo2}.


By {\bf Algorithm 2} we come up with the desired bijection. 
\begin{thm}\label{c2t_b}
There is a bijection $\phi^{(B)}_C: \mathcal{CUD}^{(B)}_n \longrightarrow \mathcal{T}^\circ_n$. 
Also, $\phi^{(B)}_C$ induces a bijection between $\mathcal{CUD}^{(B)}_{n,k}$ and $\mathcal{T}^\circ_{n,k}$ for each $1\leq k\leq n$.
\end{thm}

\begin{proof}
Given $\sigma=(a_{1,1},a_{1,2},\ldots)(a_{2,1},a_{2,2},\ldots)\ldots(a_{m,1},a_{m,2},\ldots)\in\mathcal{CUD}^{(B)}_n$. 
By applying {\bf Algorithm 2} on each cycle we obtain 
 a sequence of trees $\tau_1,\tau_2,\ldots,\tau_m$ such that $\tau_i\in\mathcal{T}^{(1)}_{\mathsf{int}(C_i)}$ with the right child of root of each $\tau_i$ being empty.  
Now we replace the empty right child of the root of $\tau_i$ by the tree $\tau_{i+1}$ for $1\leq i\leq m$.
In this way we get a tree $\tau\in\mathcal{T}^\circ_n$ and hence a bijection between $\mathcal{CUD}^{(B)}_n$ and $\mathcal{T}^\circ_n$. 
Furthermore, if $\sigma\in\mathcal{CUD}^{(B)}_{n,k}$ then $a_{m,1}=k$, hence the root of $\tau_m$ is labelled by $k$. Since the root of $\tau_m$ is the rightmost labelled node of $\tau$, we know that $\tau\in\mathcal{T}^\circ_{n,k}$. 
\end{proof}

\begin{figure}[ht]
\begin{center}
\includegraphics[scale=0.2]{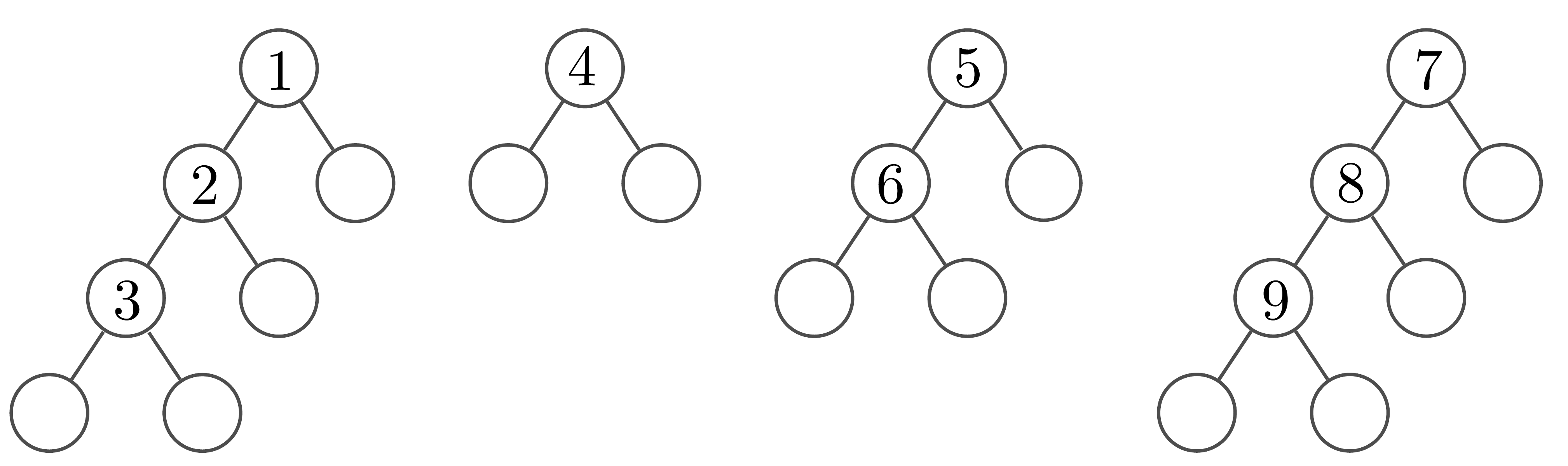}
\end{center}
\caption{\small Map each cycle to a non-plane complete increasing binary tree.}
\label{CBT1}
\end{figure}

Now we give an example in detail for the whole bijection process. Consider the permutation $\sigma=(1,\bar{3},\bar{2})(4)(5,\bar{6})(7,9,\bar{8})\in\mathcal{CUD}^{(B)}_{9,7}$. By applying {\bf Algorithm 1} on each cycle we have the bracketing of cycles
$$(1,3,2)\to(1(2(3))), \quad (4)\to(4), \quad (5,6)\to(5(6)),\quad (7,9,8)\to(7(8(9)))$$
 and therefore their corresponding non-plane trees as shown in Figure~\ref{CBT1}.

\begin{figure}[ht]
\begin{center}
\includegraphics[scale=0.2]{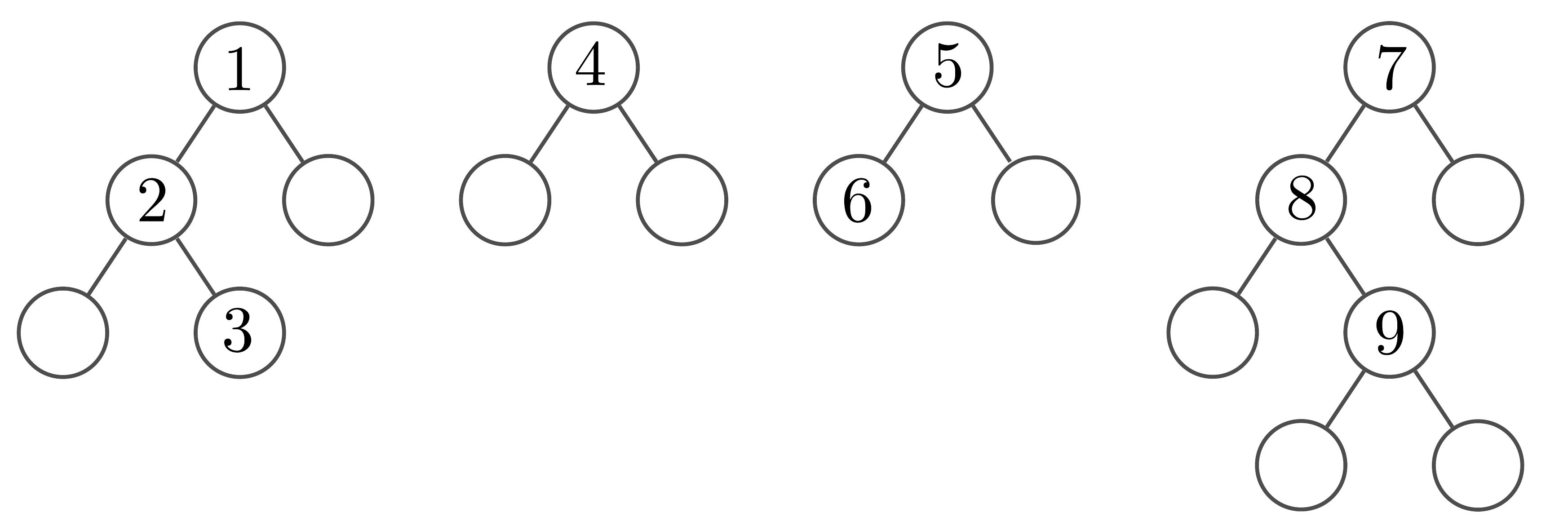}
\end{center}
\caption{\small Modify non-plane trees to plane trees $\tau_1,\tau_2,\tau_3,$ and $\tau_4$.}
\label{CBT2}
\end{figure}

Then we apply {\bf Algorithm 2} and result in the plane trees $\tau_1,\tau_2,\tau_3,$ and $\tau_4$, as shown in Figure~\ref{CBT2}.\\

\begin{figure}[ht]
\begin{center}
\includegraphics[scale=0.2]{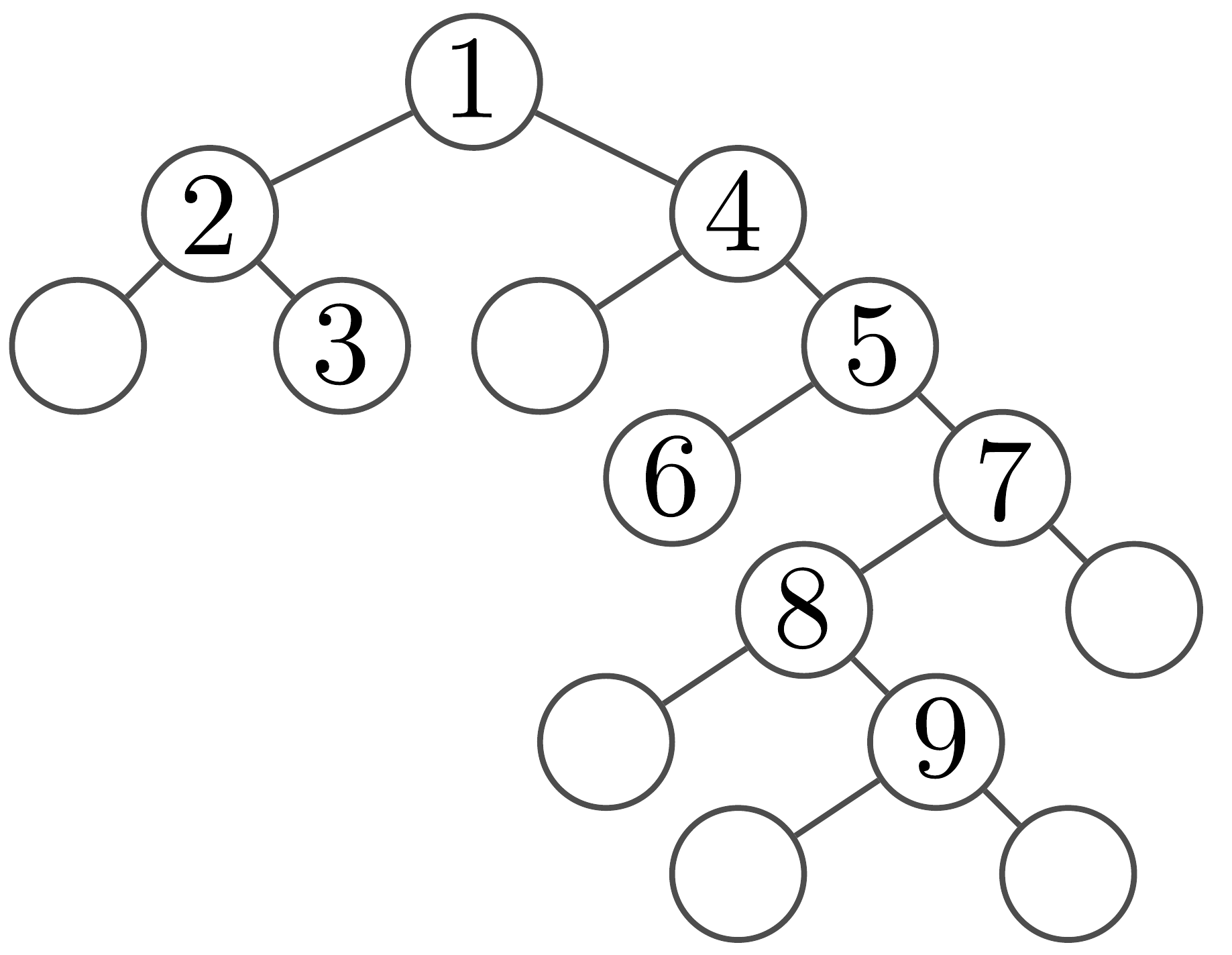}
\end{center}
\caption{\small Combine $\tau_1,\tau_2,\tau_3,$ and $\tau_4$ into a tree $\tau\in\mathcal{T}^\circ_{9,7}$.}
\label{CBT3}
\end{figure}

 Finally, using $\tau_1,\tau_2,\tau_3,$ and $\tau_4$, the tree $\tau\in\mathcal{T}^\circ_{9,7}$ can be constructed as shown in Figure~\ref{CBT3}.\\

\subsection{Type $D_n$}
The above can be modified to give a bijection in type $D_n$.

\begin{thm}\label{c2t_d}
There is a bijection $\phi^{(D)}_C: \mathcal{CUD}^{(D)}_n \longrightarrow \mathcal{T}^*_n$. 
Moreover, $\phi^{(D)}_C$ induces a bijection between $\mathcal{CUD}^{(D)}_{n,k}$ and $\mathcal{T}^*_{n,k}$
for each $1\leq k\leq n$.
\end{thm}
\begin{proof}

Let $$\sigma=(a_{1,1},a_{1,2},\ldots)(a_{2,1},a_{2,2},\ldots)\ldots(a_{m-1,1},a_{m-1,2},\ldots)(k,\bar{k})\in\mathcal{CUD}^{(D)}_{n,k}.$$
For a cycle $C_i$ of $\sigma$, do the followings.
\begin{enumerate}
    \item For $1\le i\le m-1$ : Apply {\bf Algorithm 2} on $C_i$. 
    \item For $C_m=(k,\bar{k})$ : Let $\tau_m$ be the tree of single node with label $k$.
\end{enumerate}

In this way we obtain a sequence of trees $\tau_1,\tau_2,\ldots,\tau_{m-1}$ such that $\tau_i\in\mathcal{T}^{(1)}_{\mathsf{int}(C_i)}$ for $1\leq i\leq m-1$ and 
a single node tree $\tau_m$ labeled by $k$. Since 
 for $1\leq i\leq m-1$ the right child of root of $\tau_i$ is empty, 
we replace the empty right child of the root of $\tau_i$ by the tree $\tau_{i+1}$ and get a tree $\tau\in\mathcal{T}^*_{n,k}$. This results in a bijection
between $\mathcal{CUD}^{(D)}_{n,k}$ and $\mathcal{T}^*_{n,k}$, also 
between $\mathcal{CUD}^{(D)}_n$ and $\mathcal{T}^*_n$. 
\end{proof}

\begin{figure}[ht]
\begin{center}
\includegraphics[scale=0.2]{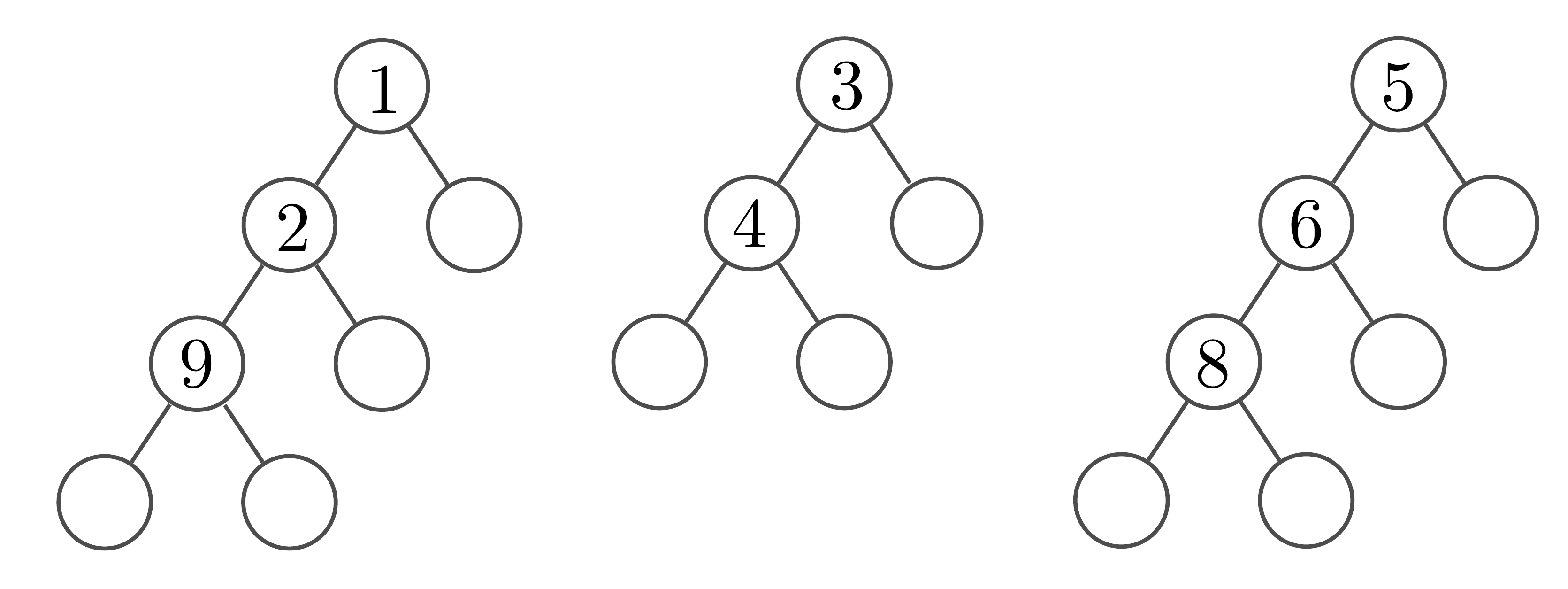}
\end{center}
\caption{\small Map each cycle to a non-plane complete increasing binary tree.}
\label{CDT1}
\end{figure}

Again we give an example in detail. Let $\sigma=(1,\bar{9},\bar{2})(3,4)(5,8,\bar{6})(7,\bar{7})\in\mathcal{CUD}^{(D)}_{9,7}$. By {\bf Algorithm 1}, the bracketing corresponding to first three cycles are 
$$(1,9,2)\to(1(2(9))),\quad  (3,4)\to(3(4)), \quad (5,8,6)\to(5(6(8)))$$
and their corresponding non-plane trees are given in Figure~\ref{CDT1}.

\begin{figure}[ht]
\begin{center}
\includegraphics[scale=0.2]{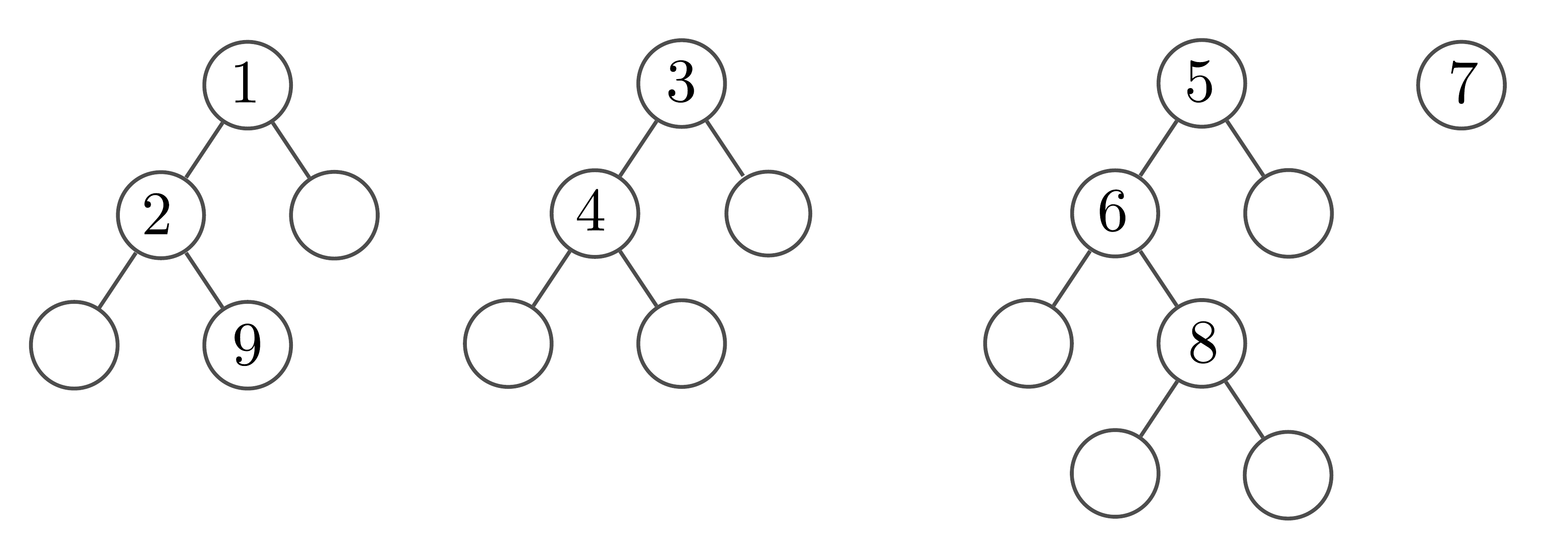}
\end{center}
\caption{\small Modify non plane trees $\tau_1,\tau_2,\tau_3$ and add a single node tree $\tau_4$. }
\label{CDT2}
\end{figure}

By applying {\bf Algorithm 2} we obtain plane trees $\tau_1,\tau_2,$ and $\tau_3$, shown in Figure~\ref{CDT2} together with the tree $\tau_4$ with a single node corresponding to the cycle $(7,\bar{7})$. 

\begin{figure}[ht]
\begin{center}
\includegraphics[scale=0.2]{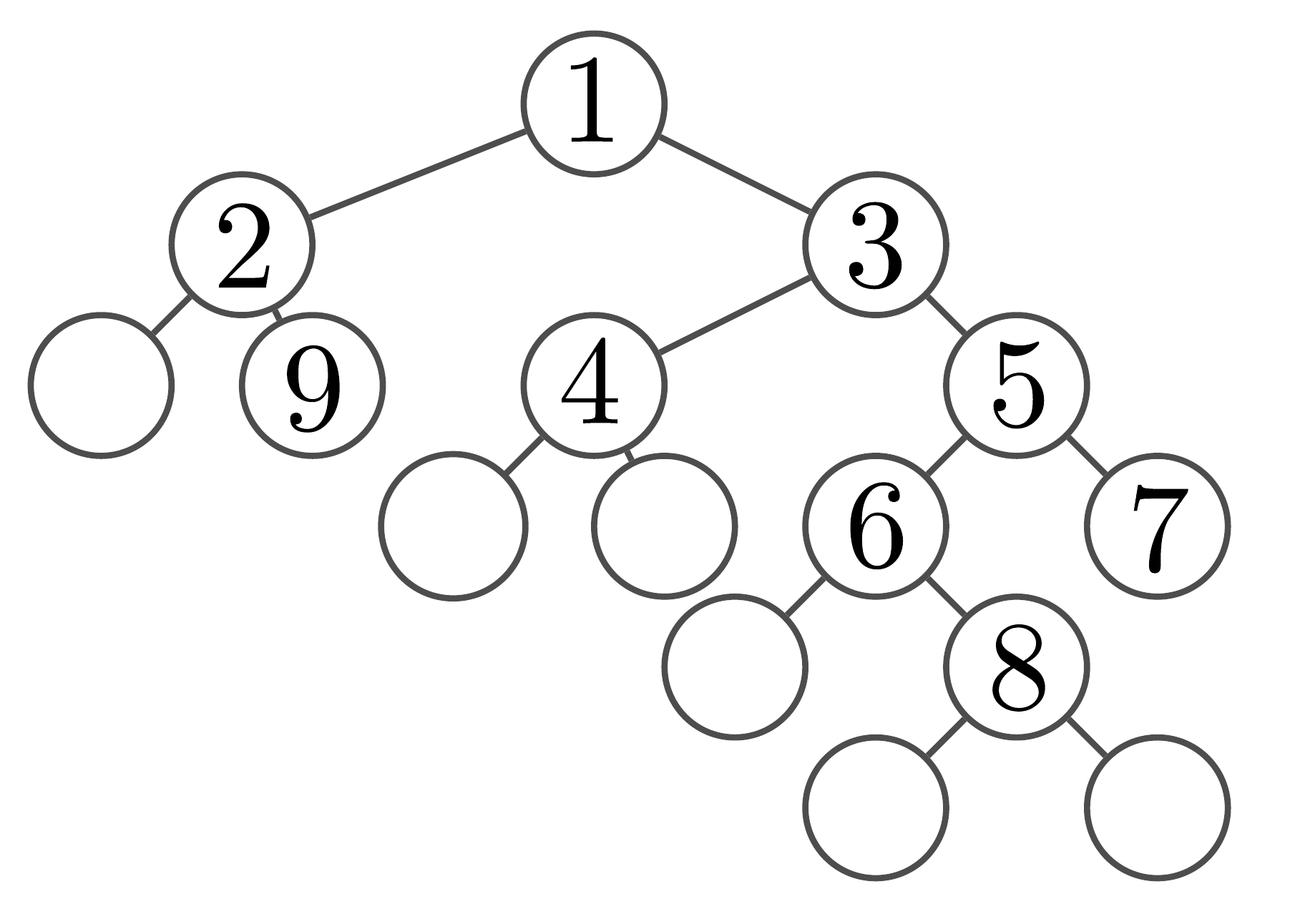}
\end{center}
\caption{\small Combine $\tau_1,\tau_2,\tau_3,$ and $\tau_4$ into a tree $\tau\in\mathcal{T}^*_{9,7}$.}
\label{CDT3}
\end{figure}

Finally, the tree $\tau\in\mathcal{T}^*_{9,7}$ is constructed from $\tau_1,\tau_2,\tau_3,$ and $\tau_4$, 
as shown in Figure~\ref{CDT3}.\\

The following bonus is straightforward  from Theorem~\ref{c2t_b} and Theorem~\ref{c2t_d}
and we omit the proof.
\begin{cor}
In both bijections $\phi^{(B)}_C:\sigma\mapsto\tau$ and $\phi^{(D)}_C:\sigma\mapsto\tau$, the labels on the rightmost path of $\tau$ is exactly the set of minimum entries of cycles of $|\sigma|$.
\end{cor}

\section{Complete increasing binary Trees and Valley signed permutations}
In this section, we give the bijections 
 $$\phi^{(B)}_V:\mathcal{VS}^{(B)}_n\to\mathcal{T}^\circ_n \quad \text{and} \quad \phi^{(D)}_V:\mathcal{VS}^{(D)}_n\to\mathcal{T}^*_n$$
between valley signed permutations and complete increasing binary trees. 

\subsection{Type $B_n$}

We first look at the type $B_n$ case. Firstly, note that for $\sigma=\sigma_1\sigma_2\ldots\sigma_n\in\mathfrak{S}_n$
in its window notation, if $\{\sigma_{v_1},\sigma_{v_2},\ldots,\sigma_{v_t}\}$ are the set of valleys, then a valley signed permutation $\sigma'\in\mathcal{VS}^{(B)}_n$ with $|\sigma'|=\sigma$ is uniquely determined from the choices of the signs of $\sigma'_{v_1+1},\sigma'_{v_2+1},\ldots,\sigma'_{v_t+1}$. For example, there are $4$
type $B_n$ valley signed permutations
$51324, 51\bar{3}24, 5132\bar{4}, 51\bar{3}2\bar{4} \in \mathcal{VS}^{(B)}_5$ such that
$|\sigma'|=51324$ which has two valleys (at $1$ and $2$).

The next algorithm, which sends an ordinary permutation to a complete increasing binary plane tree with empty leaves, is well-known (\cite{EC1}, p.44).\\ 

\noindent \textbf{Algorithm 3}\\
\noindent \textbf{Input :} A sequence $\sigma=\sigma_1,\sigma_2,\ldots,\sigma_n$ of distinct positive integers.

\medskip
Let $\sigma_i=\min(\sigma)$ and $\sigma_L=\sigma_1,\ldots,\sigma_{i-1}$ and $\sigma_R=\sigma_{i+1},\ldots,\sigma_n$.
Construct a tree recursively by setting $\sigma_i$ as the root of $\mathcal{T}$ and $\mathcal{T}(\sigma_L)$ (resp. $\mathcal{T}(\sigma_R$)) be its right and left subtree (with the convention that $T(\emptyset)$ is an empty node).
Note that to our purpose we set $\mathcal{T}(\sigma_L)$ to be the `right' subtree.

\medskip
\noindent \textbf{Output :} A complete increasing binary plane tree with empty leaves.\\

\begin{figure}[ht]
\begin{center}
\includegraphics[scale=0.2]{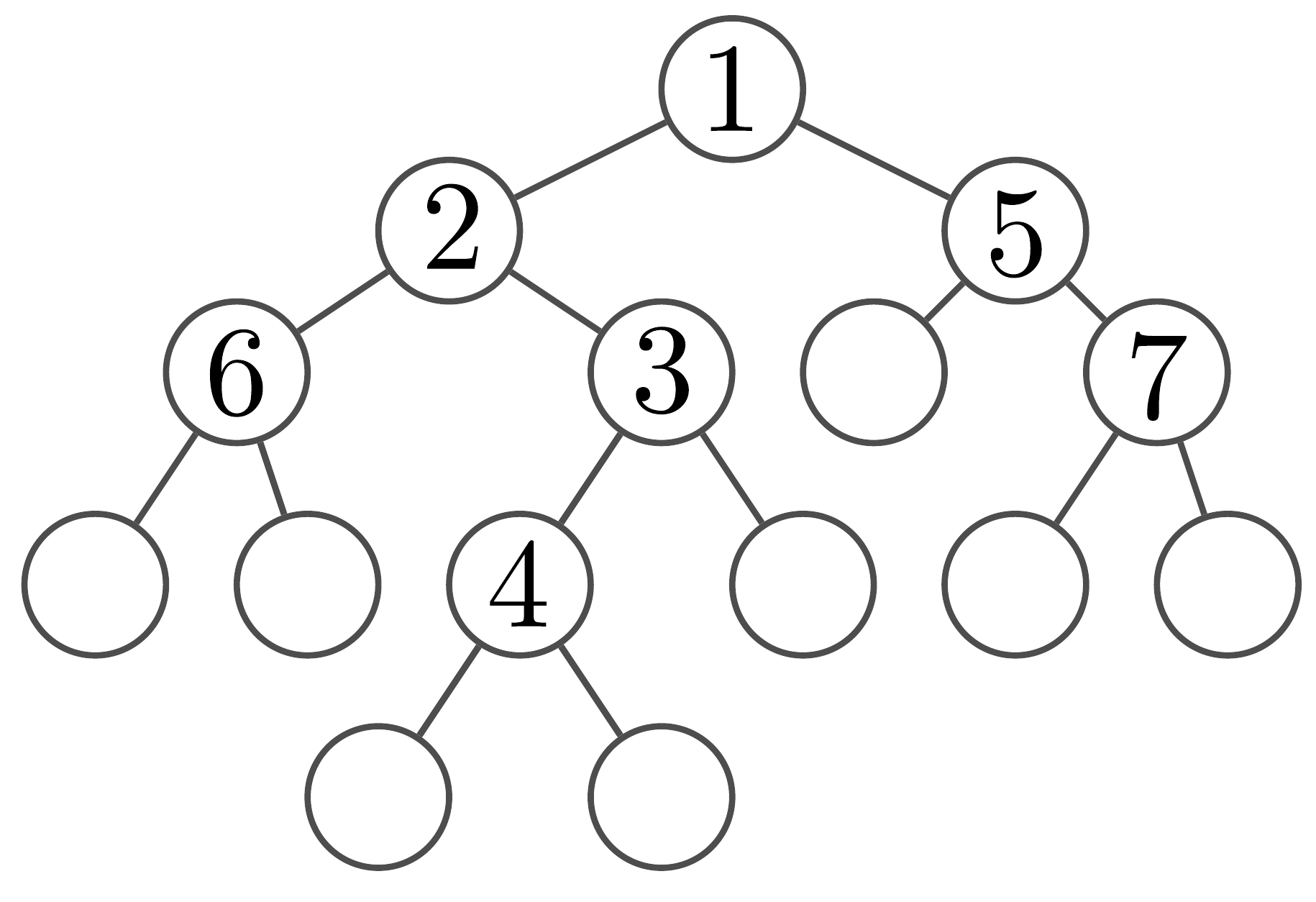}
\end{center}
\caption{\small The tree obtained from $7513426$ by {\bf Algorithm 3}.}
\label{algo3}
\end{figure}

For example, the tree obtained from $7513426$ is shown in Figure~\ref{algo3}.

From {\bf Algorithm 3} we have an easy observation.
\begin{lem}\label{peakleaf}
For $2\leq i\leq n$, the node with label $\sigma_i$ has two empty children in $\mathcal{T}(\sigma)$ 
if and only if $\sigma_i$ is a peak in $\sigma$. 
\end{lem}
\begin{proof}
For $2\leq i\leq n$, if $\sigma_i$ is a peak, then it is not the minimum entry of any consecutive subsequence of $\sigma$ unless the subsequence only contains $\sigma_i$ itself. This means that when $\mathcal{T}(\sigma)$ is constructed both of the left and right subtrees of $\sigma_i$ are empty nodes. If $\sigma_i$ is not a peak, then at least one of its subtree is nonempty.
\end{proof}

\begin{thm}\label{v2t_b}
There is a bijection $\phi^{(B)}_V: \mathcal{VS}^{(B)}_n \longrightarrow \mathcal{T}^\circ_n$.
The map $\phi^{(B)}_V$ also induces a bijection between $\mathcal{VS}^{(B)}_{n,k}$ and $\mathcal{T}^\circ_{n,k}$ for each
$1\leq k\leq n$.
\end{thm}

\begin{proof}
Given $\sigma\in\mathcal{VS}^{(B)}_n$, we first construct $\mathcal{T}(|\sigma|)$ by {\bf Algorithm 3}. 
If $|\sigma_{v_1}|,|\sigma_{v_2}|,\ldots,|\sigma_{v_t}|$ are the valleys of $|\sigma|$, then we have
$$\{i:\sigma_i<0\}\subseteq\{v_1+1,v_2+1,\ldots,v_t+1\}.$$ 
There must exist a peak $|\sigma_{p_i}|$ in $|\sigma|$ such that $v_i<p_i<v_{i+1}$ for each $i=1,2,\ldots,t-1$, and a peak $|\sigma_{p_t}|$ with $p_t>v_t$. Now, for $1\leq i\leq t$, if $\sigma_{v_{i}+1}<0$, we remove the two empty children of the node with label $|\sigma_{p_i}|$ in $\mathcal{T}(\sigma)$. Note that these operations are well defined by Lemma~\ref{peakleaf}.
We result in a tree $\tau\in\mathcal{T}^\circ_n$ and set $\phi^{(B)}_V(\sigma)=\tau$. 

Furthermore, since $\sigma_1$ is in the left subsequence during the recursive construction of $\mathcal{T}(\sigma)$, the node with label $\sigma_1$ must be the rightmost labelled node. Hence the map $\phi^{(B)}_V$ induces a bijection between $\mathcal{VS}^{(B)}_{n,k}$ and $\mathcal{T}^\circ_{n,k}$.
\end{proof}

\begin{figure}[ht]
\begin{center}
\includegraphics[scale=0.2]{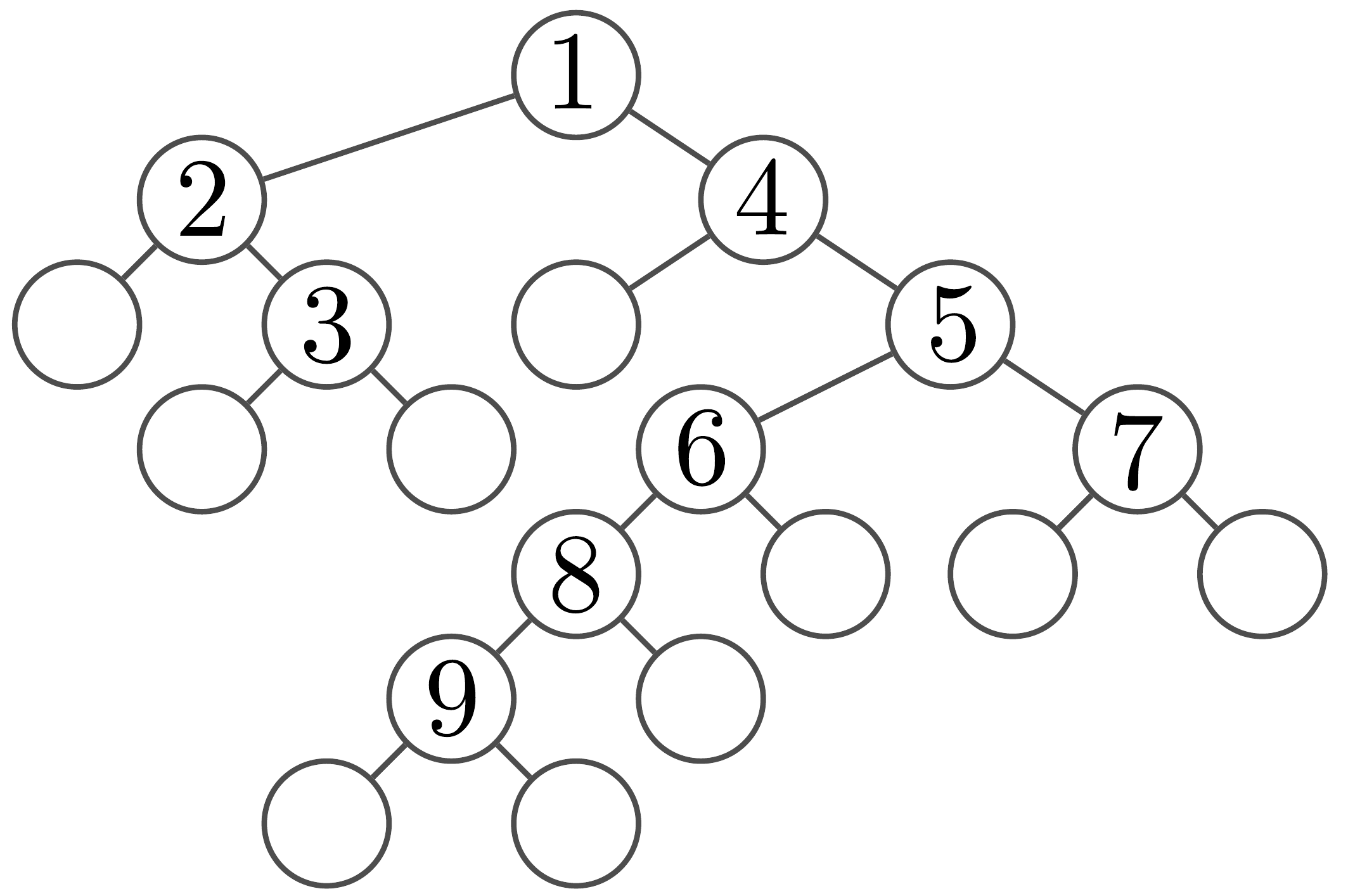}
\qquad
\includegraphics[scale=0.2]{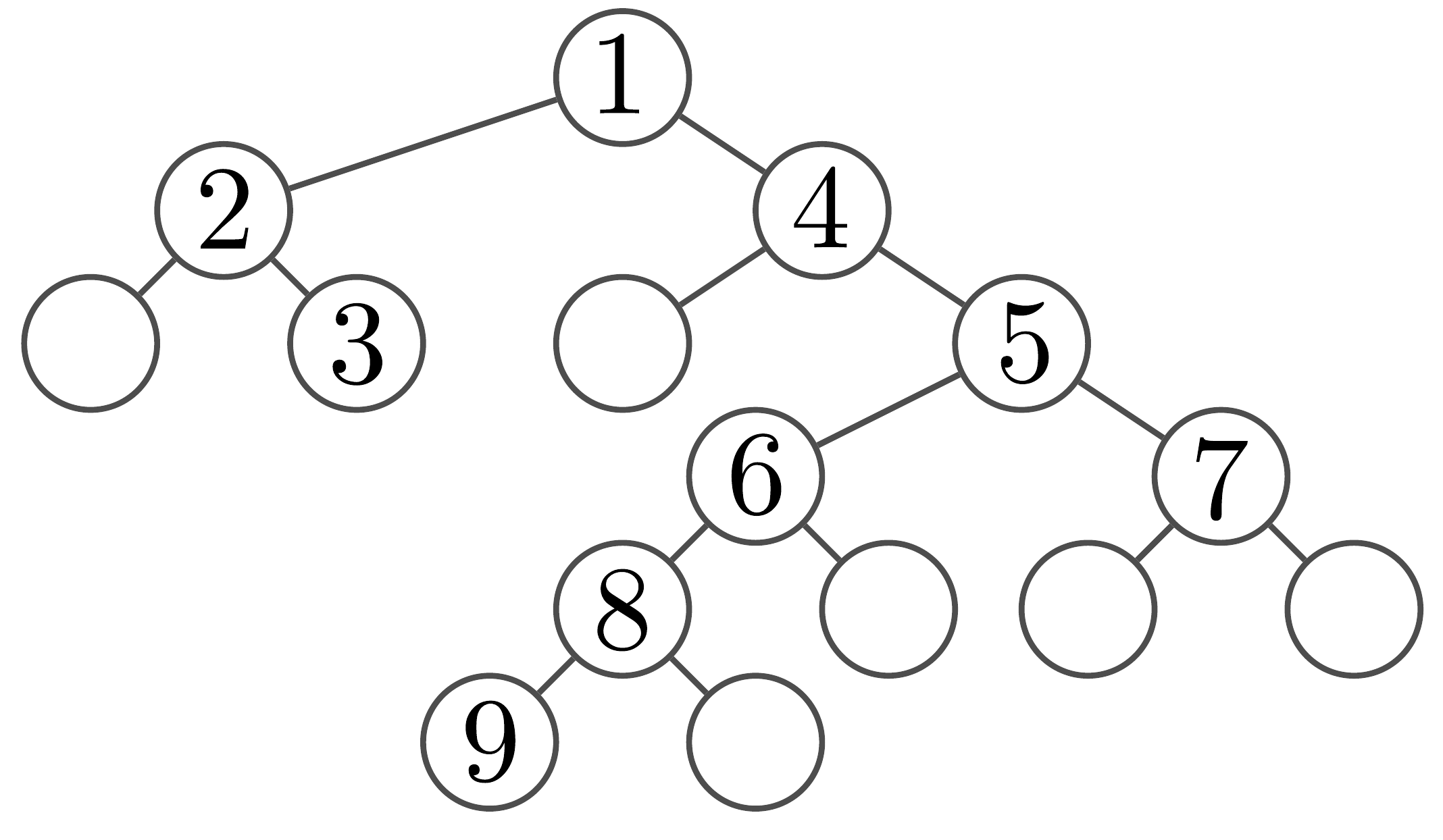}
\end{center}
\caption{\small The process of bijection from $75\bar{6}8941\bar{3}2\in\mathcal{VS}^{(B)}_{9,7}$ to a tree 
in $\mathcal{T}^\circ_{9,7}$.}
\label{VBT}
\end{figure}

For example, let $\sigma=75\bar{6}8941\bar{3}2\in\mathcal{VS}^{(B)}_{9,7}$. By {\bf Algorithm 3}, $\mathcal{T}(|\sigma|)$ is constructed as in the left of Figure~\ref{VBT}.  

Since the peaks in $|\sigma|$ corresponding to negative entries ($\bar{6}$ and $\bar{3}$) are $9$ and $3$, we remove the empty children of these two nodes from $\mathcal{T}(|\sigma|)$ and obtain the tree $\tau\in\mathcal{T}^\circ_{9,7}$ in the right of Figure~\ref{VBT}.

\subsection{type $D_n$}

Similarly to the type $B_n$ case, if $\{\sigma_{v_1},\sigma_{v_2},\ldots,\sigma_{v_t}\}$ are the set of valleys of $\sigma\in \mathfrak{S}_n$, then a valley signed permutation $\sigma'\in\mathcal{VS}^{(D)}_n$ with $|\sigma'|=\sigma$ is also uniquely determined from the choices of the signs of $\sigma'_{v_1+1},\sigma'_{v_2+1},\ldots,\sigma'_{v_t+1}$. For example,
There are four $\sigma'\in \mathcal{VS}^{(D)}_5$ with 
$|\sigma'|=51324$, namely $\bar{5}1324, \bar{5}1\bar{3}24, \bar{5}132\bar{4}$, and $\bar{5}1\bar{3}2\bar{4}$.

\begin{thm}\label{v2t_d}
There is a bijection $\phi^{(D)}_V:\mathcal{VS}^{(D)}_n \longrightarrow \mathcal{T}^*_n$. The map $\phi^{(D)}_C$ also induces a bijection between $\mathcal{VS}^{(D)}_{n,k}$ and $\mathcal{T}^*_{n,k}$ for each $1\leq k\leq n$..
\end{thm}

\begin{proof}
Given $\sigma\in\mathcal{VS}^{(D)}_n$, we first construct $\mathcal{T}(|\sigma|)$ by {\bf Algorithm 3}. 

Firstly, from the definition of $\mathcal{VS}^{(D)}_n$, the rightmost labelled node of $\mathcal{T}(\sigma)$ is labelled by $|\sigma_1|$ with two empty children and we remove them. 
Let $|\sigma_{v_1}|,|\sigma_{v_2}|,\ldots,|\sigma_{v_t}|$ be the valleys of $|\sigma|$. We have
$$\{i:i>1\mbox{ and }\sigma_i<0\}\subseteq\{v_1+1,v_2+1,\ldots,v_t+1\}.$$
There must exist a peak $|\sigma_{p_i}|$ in $|\sigma|$ such that $v_i<p_i<v_{i+1}$ for each $i=1,2,\ldots,t-1$, and a peak $|\sigma_{p_t}|$ with $p_t>v_t$. 
Now for $1\leq i\leq t$, if $\sigma_{v_{i}+1}<0$, then remove the two empty children of the node with label $|\sigma_{p_i}|$ in $\mathcal{T}(\sigma)$. These operations are well defined by Lemma~\ref{peakleaf}.
We obtain the resulting tree $\tau\in\mathcal{T}^*_n$ and set $\phi^{(D)}_V(\sigma)=\tau$.

Furthermore, since the rightmost labelled node of $\tau$ is labelled by $|\sigma_1|$, the map $\phi^{(D)}_V$ induces a bijection between $\mathcal{VS}^{(D)}_{n,k}$ and $\mathcal{T}^*_{n,k}$.
\end{proof}

\begin{figure}[ht]
\begin{center}
\includegraphics[scale=0.2]{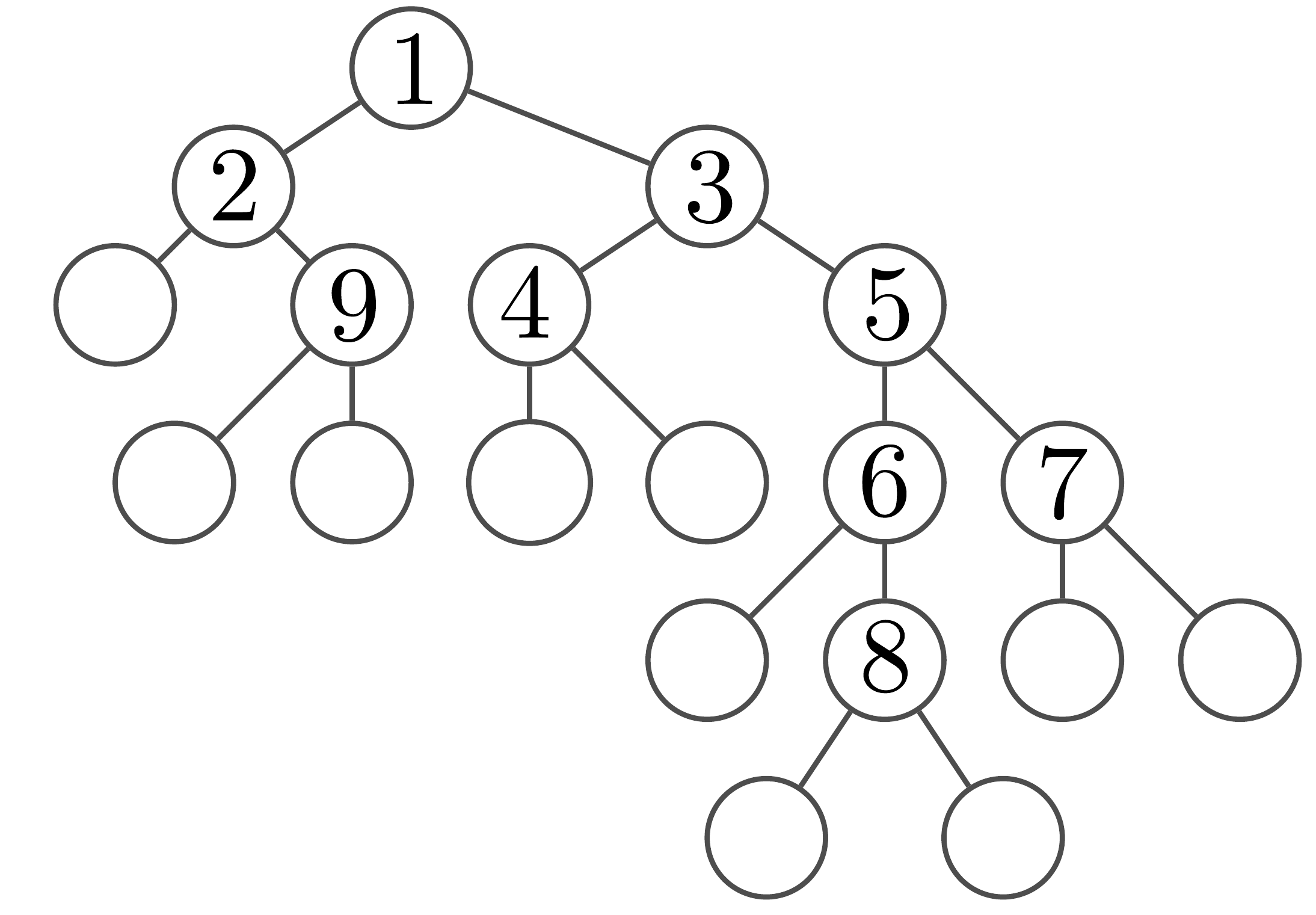}
\qquad
\includegraphics[scale=0.2]{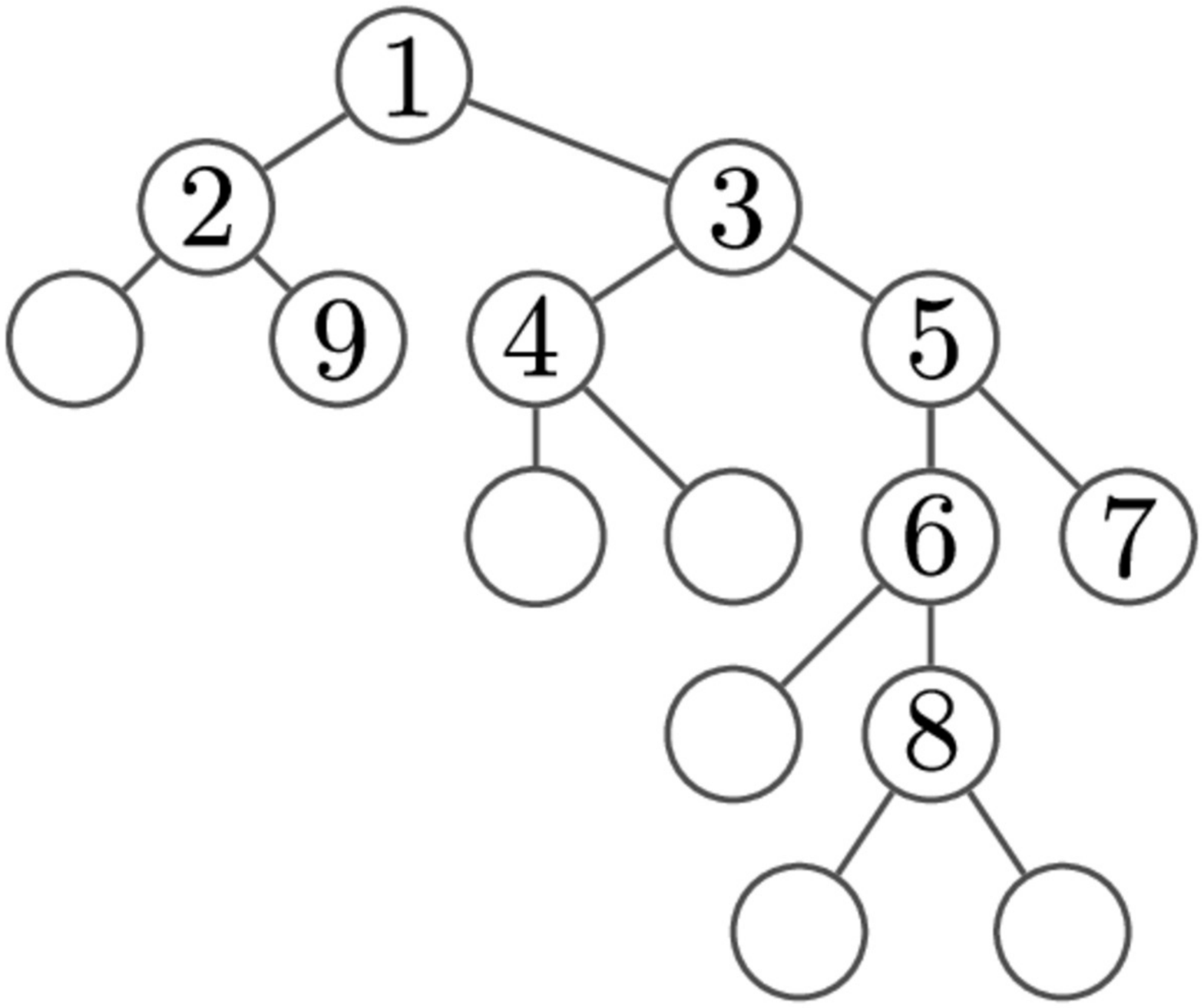}
\end{center}
\caption{\small The process of bijection from $\bar{7}586341\bar{9}2\in\mathcal{VS}^{(D)}_{9,7}$ to a tree 
in $\mathcal{T}^*_{9,7}$.}
\label{VDT}
\end{figure}

For example, let $\sigma=\bar{7}586341\bar{9}2\in\mathcal{VS}^{(D)}_{9,7}$. The tree 
$\mathcal{T}(|\sigma|)$ is shown in the left of Figure~\ref{VDT}. First we remove the two empty children of the node $7$.
The negative entries $\bar{9}$ is corresponding to the peak $9$ in $|\sigma|$, so we remove its empty children from 
$\mathcal{T}(|\sigma|)$ and obtain the tree $\tau\in\mathcal{T}^*_{9,7}$ in the right of Figure~\ref{VDT}.

For a permutation $\sigma=\sigma_1\sigma_2\ldots\sigma_n$, an entry $\sigma_i$ is a {\it left-to-right minimum} if $\sigma_i=\min\{\sigma_j:1\leq j\leq i\}$.
We have a bonus corollary. 

\begin{cor}
In both bijections $\phi^{(B)}_V$ and $\phi^{(D)}_V$, the set of the labels on the rightmost path of $\tau$ is the set of left-to-right minimums of $|\sigma|$.
\end{cor}
\begin{proof}
Observe that in $\phi^{(B)}_V$ or $\phi^{(D)}_V$,  $|\sigma_i|$ is a label on the rightmost path if and only if $|\sigma_i|$ is on the left of the minimum entry at each stage of {\bf Algorithm 3} until itself becomes the minimum. Hence the entries on the left of $|\sigma_i|$ are all with absolute value greater than $|\sigma_i|$. That is, $|\sigma_i|$ is a left-to-right minimum.
\end{proof}


\section{Proof of Main Theorem III}
We turn to flip equivalence classes. For completeness' sake we 
give a proof of Knuth's result in terms of increasing 1-2 trees. An \emph{increasing 1-2 tree} on the nodes $[n]$ is a  (non-plane) rooted tree so that each non-leaf nodes has $1$ or $2$ children and labels along the unique path from  root to any leaf are increasing. It is well known that they are also counted by Euler numbers~\cite{Stanley_10,EC1}.\\

\textit{Proof of Proposition~\ref{flipSn}}: It suffices to maps an equivalent class to an increasing 1-2 tree.
For $\sigma\in\mathfrak{S}_n$, let $\hat{\sigma}=\min(\sigma)$ and write $\sigma=\sigma_L\hat{\sigma}\sigma_R$. Construct an increasing 1-2 tree $\tau(\sigma)$ recursively as follows:
\begin{enumerate}
    \item Let $\hat{\sigma}$ be  the root of $\tau(\sigma)$.
    \item Let $\tau(\sigma_L)$ and $\tau(\sigma_R)$ be the (non-ordered) subtrees of the root $\hat{\sigma}$.
\end{enumerate}
It can be checked directly that flip equivalent permutations map to the same tree and the result follows. 
\hfill $\Box$

\subsection{Bijection between flip equivalence classes and complete binary trees}
The main result of this subsection is the next bijection.
\begin{thm}\label{f2t}
 There is a bijection $$\phi_F:(\mathcal{FL}^{(B)}_n\cup\mathcal{FL}^{(D)}_n)\to(\mathcal{T}^\circ_n\cup\mathcal{T}^*_n).$$
Moreover, $\phi_F$ induces a bijection between $\mathcal{FL}^{(B)}_{n,k}$ ($\mathcal{FL}^{(D)}_{n,k}$, respectively) and $\mathcal{T}^\circ_{n,k}$ ($\mathcal{T}^*_{n,k}$, respectively) for $1\leq k\leq n$.
\end{thm}

\begin{proof}
We shall proceed as follows. For a given signed permutation $\sigma$ we first correspond it to a tree $\tau(\sigma)\in\mathcal{T}^\circ_n\cup\mathcal{T}^*_n$. Then we show that if $\sigma'\in[\sigma]$, then $\tau(\sigma')=\tau(\sigma)$ and therefore the bijection $\phi_F([\sigma]):=\tau(\sigma)$ is well-defined.\\

\noindent \textbf{Algorithm 4}\\
\noindent \textbf{Input :} A signed permutation $\sigma$. 

We construct the tree $\tau(\sigma)\in\mathcal{T}^\circ_n\cup\mathcal{T}^*_n$ recursively.
If $\sigma$ is empty, let $\tau(\sigma)$ be an empty node.
If $\sigma$ is not empty, write $\sigma=\sigma_L,\sigma_i,\sigma_R$, where $|\sigma_i|=\min(|\sigma|)$. Let the label of the root of $\tau(\sigma)$ be $|\sigma_i|$.
    \begin{enumerate}
        \item If both of $\sigma_L,\sigma_R$ are empty and $\sigma_i<0$, then let the node with label $|\sigma_i|$ be a leaf.
        \item If both of $\sigma_L,\sigma_R$ are empty and $\sigma_i>0$, then let the node with label $|\sigma_i|$ be with two empty leaves.
        \item If $\sigma_L$ and $\sigma_R$ are not both empty, consider the two values $\min(|\sigma_L|)$ and $\min(|\sigma_R|)$, where $\min(\emptyset):=+\infty$.
        \begin{itemize}
            \item[---] If $\sigma_i>0$ and $\min(|\sigma_L|)<\min(|\sigma_R|)$, then let $\tau(\sigma_L)$ ($\tau(\sigma_R)$, respectively) be the left (right, respectively) subtree of the node with label $|\sigma_i|$.
            \item[---] If $\sigma_i>0$ and $\min(|\sigma_L|)>\min(|\sigma_R|)$, then let $\tau(\sigma_R)$ ($\tau(\sigma_L)$, respectively) be the left (right, respectively) subtree of the node with label $|\sigma_i|$.
            \item[---] If $\sigma_i<0$ and $\min(|\sigma_L|)<\min(|\sigma_R|)$, then let $\tau(\sigma_R)$ ($\tau(\sigma_L)$, respectively) be the left (right, respectively) subtree of the node with label $|\sigma_i|$.
            \item[---] If $\sigma_i<0$ and $\min(|\sigma_L|)>\min(|\sigma_R|)$, then let $\tau(\sigma_L)$ ($\tau(\sigma_R)$, respectively) be the left (right, respectively) subtree of the node with label $|\sigma_i|$.
        \end{itemize}
    
    \end{enumerate}
\noindent \textbf{Output :} A complete increasing binary tree $\tau(\sigma)$.\\

It is clear that $\tau(\sigma)\in\mathcal{T}_n$. By the definition of flip equivalence of signed permutations, it is clearly that  the bijection $\phi_F$ defined by $\phi_F([\sigma]):=\tau(\sigma)$ is well-defined. Moreover, the node with label $|\mathsf{smax}([\sigma])|$ is the rightmost label of the tree $\phi_F([\sigma])$. That is, for $1\leq k\leq n$, the map $\phi_V$ induces a bijection between $\mathcal{FL}^{(B)}_{n,k}$ ($\mathcal{FL}^{(D)}_{n,k}$, respectively) and $\mathcal{T}^\circ_{n,k}$ ($\mathcal{T}^*_{n,k}$, respectively).
\end{proof}

\begin{figure}[ht]
\begin{center}
\includegraphics[scale=0.2]{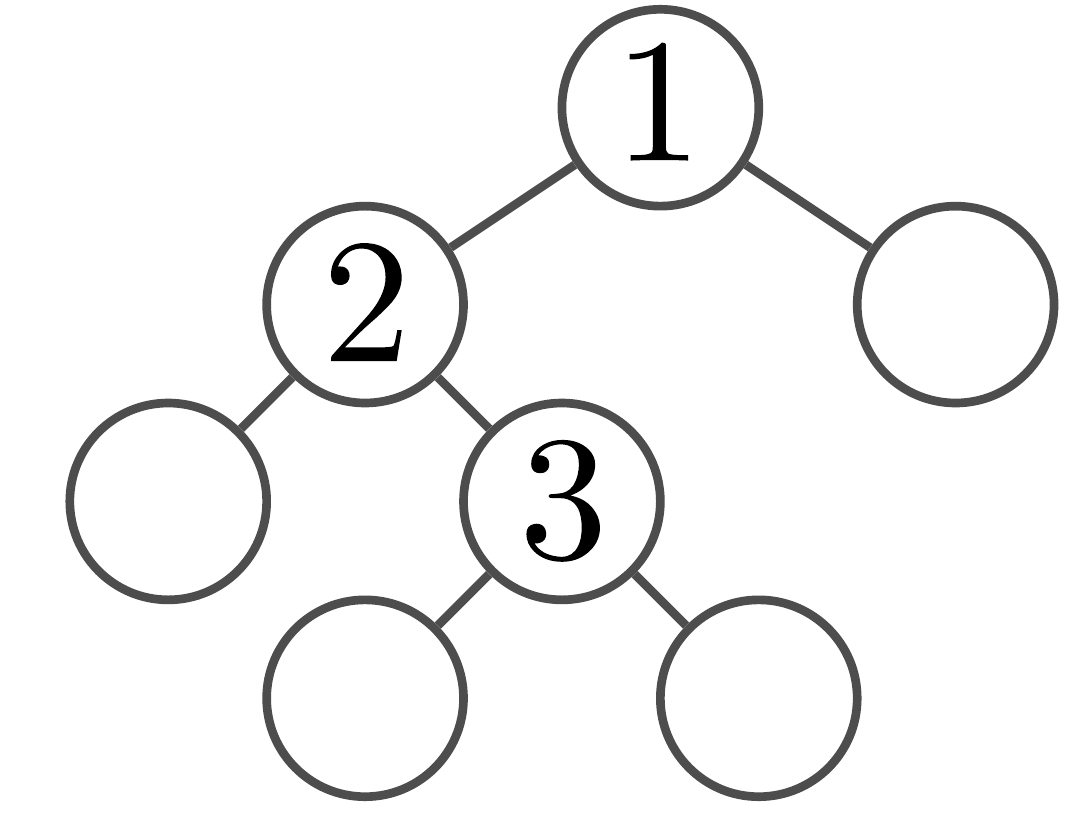}
\qquad
\includegraphics[scale=0.2]{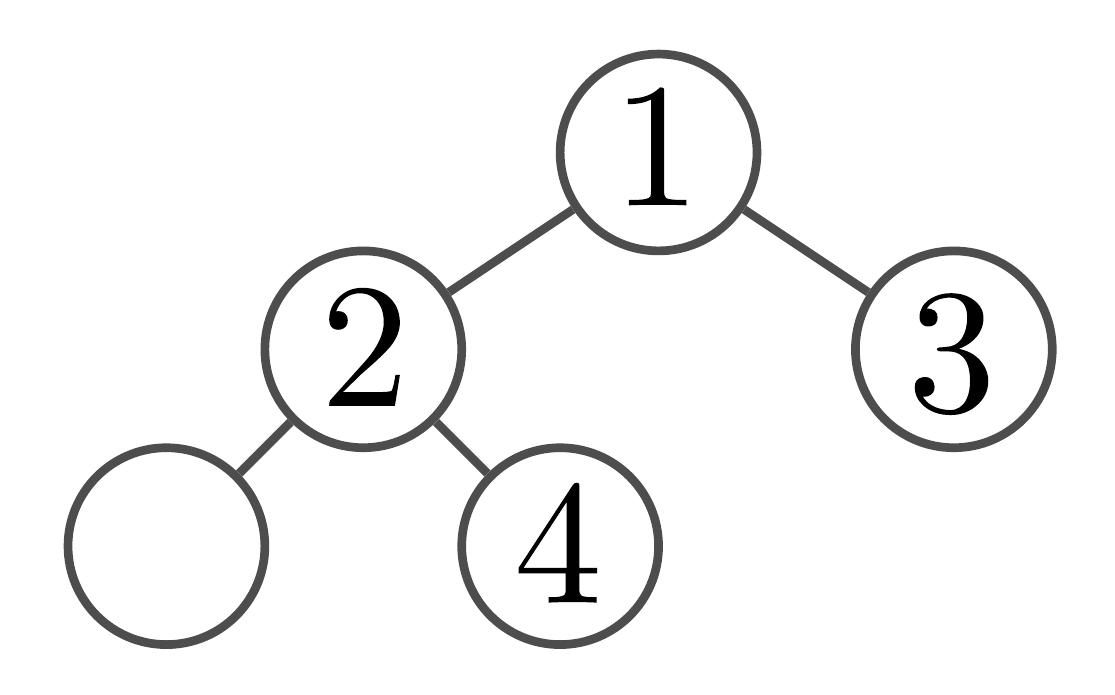}
\end{center}
\caption{\small The trees correspond to $[1\bar{2}3]$ and $[\bar{2}\bar{4}13]$.}
\label{FT}
\end{figure}

For example, the equivalence class $[1\bar{2}3]\in\mathcal{FL}^{(B)}_{3,1}$ contains $4$ signed permutations $1\bar{2}3,3\bar{2}1,\bar{2}31,$ and $13\bar{2}$. Apply bijection $\phi_F$ on each of the $4$ signed permutations we get the same tree $\tau\in\mathcal{T}^\circ_{3,1}$, as shown in the left of Figure~\ref{FT}. 
On the other hand, the equivalence class
$[\bar{2}\bar{4}13]\in \mathcal{FL}^{(D)}_{4,3}$ contains $4$ permutations
$\bar{2}\bar{4}13, \bar{3}1\bar{2}\bar{4}, \bar{3}1\bar{4}\bar{2}$ and $\bar{4}\bar{2}1\bar{3}$ and the tree we obtained is in the right of Figure~\ref{FT}. 

\subsection{Proof of Theorem 2.4}
To prove Main Theorem III, we need a final piece.

\begin{lem}\label{emp2nva}
For $[\sigma]\in\mathcal{FL}^{(B)}_n\cup\mathcal{FL}^{(D)}_n$, we have
$$\mathsf{emp}(\phi_F([\sigma]))=n-2\cdot\mathsf{spk}([\sigma])+1.$$
\end{lem}
\begin{proof}
Observe that $\mathsf{spk}([\sigma])$ is equal to the number of labelled leaves of $\phi_F([\sigma])$. Hence the number of non-leaf nodes of $\phi_F([\sigma])$ is $n-\mathsf{spk}(\sigma)$. Moreover, the total number of leaves of a complete binary tree is one plus the number of its non-leaf nodes. Therefore, we have $$\mathsf{emp}(\phi_F([\sigma]))+\mathsf{spk}([\sigma])=n-\mathsf{spk}([\sigma])+1$$ which is our goal.   
\end{proof}
\begin{proof}[Proof of Theorem~\ref{flpoly}.]
By Theorem~\ref{treeV}, Theorem~\ref{f2t} and Lemma~\ref{emp2nva} the proof is completed.  
\end{proof}


\section{Concluding remarks}
In this paper, we give types $B_n$ and type $D_n$ definitions of cycle-up-down permutations and
flip equivalence classes on signed permutations. We provide new refined Arnold families on them as well as on the 
valley signed permutations. Table~\ref{summary} summarizes our contribution.
 
\begin{center}
\begin{table}[ht]

\begin{tabular}{l|ll|ll|ll|}
\cline{2-7}
                       & \multicolumn{2}{l|}{ \qquad  Cycle-up-down perm.}    & \multicolumn{2}{l|}{  \qquad Valley-signed perm.}    & \multicolumn{2}{l|}{  \qquad Flip equiv. classes}    \\ \cline{2-7} 
                       & \multicolumn{1}{l|}{Definition} & \qquad $V_{n,k}$  & \multicolumn{1}{l|}{Definition} &  \qquad   $V_{n,k}$ & \multicolumn{1}{l|}{Definition} &  \qquad  $V_{n,k}$ \\ \hline

\multicolumn{1}{|l|}{Type $A_{n-1}$} & \multicolumn{1}{l|}{ \qquad \cite{Elizalde_Deutsch_11}} & { \qquad ---}  & \multicolumn{1}{l|}{ \quad ---} & { \qquad ---}  & \multicolumn{1}{l|}{ \quad \cite{Knuth_09}} & { \qquad ---}  \\ 

\multicolumn{1}{|l|}{Type $B_n$} & \multicolumn{1}{l|}{Section \ref{CUDBD}} & Theorem~\ref{cudpoly}  & \multicolumn{1}{l|}{ \quad \cite{JosuatVerges_Novilli_Thibon_12}} & Theorem~\ref{vspoly}  & \multicolumn{1}{l|}{Section \ref{KNUBD}} &  Theorem~\ref{flpoly} \\ 

\multicolumn{1}{|l|}{Type $D_n$} & \multicolumn{1}{l|}{Section \ref{CUDBD}} & Theorem~\ref{cudpoly}  & \multicolumn{1}{l|}{ \quad \cite{JosuatVerges_Novilli_Thibon_12}} & Theorem~\ref{vspoly}  & \multicolumn{1}{l|}{Section \ref{KNUBD}} &  Theorem~\ref{flpoly} \\ \hline
\end{tabular}
\smallskip
\caption{Summary of results.}
~\label{summary}
\end{table}
\end{center}

We also give bijections between these structures with complete increasing binary trees. Together with those bijections in~\cite{Eu_Fu_23} between complete increasing binary trees and other structures (e.g. signed increasing 1-2 trees, signed Andr\'{e} permutations, and signed Simsun permutations of type I and type II), we  have bijections between all these Arnold families. 

However, in the light of the fact that
the set of even-signed permutations is just a subset of signed permutations, i.e., $\mathfrak{D}_n \subset \mathfrak{B}_n$, it is not so satisfying that in our structures as well as those presented in~\cite{Eu_Fu_23} realizing the refined Arnold families, the type $D_n$ structures are usually not subsets of the corresponding type $B_n$ structures. 
This hints that there could exist totally different constructions of all these structures compatible with $\mathfrak{D}_n \subset \mathfrak{B}_n$. It will be of great interest if one can find them.

Euler numbers are ubiquitous in combinatorics. Except for what we and~\cite{Eu_Fu_23} have investigated, there are still many interesting combinatorial structures counted by them. For examples, the number of sequences $(a_1, \dots a_{n-1}$) with $0 \le a_i < i$, such that no three terms are equal~\cite{Corteel_Martinez_Savage_Weselcouch_16}, the number of self-dual rooted edge-labeled trees with $n$ vertices~\cite{Apostolakis_18}, the number of
linear extensions of the ``zig-zag" poset~\cite{Stanley_10,EC1}, the number of total cyclic orders on $\{0,1,\dots ,n\}$ under certain conditions~\cite{Ramassamy_17}, and so on. It would also be interesting to find their types $B_n$ and $D_n$ counterpart definitions and corresponding statistics such that together they realize the refined Arnold families. 

Another common generalization of $\mathfrak{S}_n$ and $\mathfrak{B}_n$ is the wreath product $\mathfrak{S}_n\wr C_n$, the $r$-colored permutations. The alternating permutations have been generalized along this route~\cite{JosuatVerges_Novilli_Thibon_12}. It is also a natural problem to see if all the structures considered in~\cite{Eu_Fu_23} and this paper have their colored generalizations. We leave these problems to the interested readers.



\end{document}